\documentclass[10pt]{article}
\usepackage[utf8]{inputenc}
\usepackage{amsmath,amsthm,amsfonts}
\usepackage{graphicx}
\usepackage{color}
\usepackage{multirow}

\graphicspath{ {./figs/}}

\newcommand{\const}{\mathop{\rm const}\nolimits}

\begin{document}

\title{Learning macroscopic parameters in nonlinear multiscale simulations using nonlocal multicontinua upscaling techniques}

\author{
Maria Vasilyeva\thanks{Institute for Scientific Computation, Texas A\&M University, College Station, TX, USA  \& Multiscale model reduction laboratory, North-Eastern Federal University, Yakutsk, Russia. Email: {\tt vasilyevadotmdotv@gmail.com}.}
\and
Wing T. Leung\thanks{ICES, University of Texas, Austin, TX, USA (\texttt{wleung@ices.utexas.edu})}
\and
 Eric T. Chung\thanks{Department of Mathematics, The Chinese University of Hong Kong, Shatin, New Territories, Hong Kong SAR, China (\texttt{tschung@math.cuhk.edu.hk}) }
\and
Yalchin Efendiev\thanks{Department of Mathematics \& Institute for Scientific Computation (ISC), Texas A\&M University,
College Station, TX, USA (\texttt{efendiev@math.tamu.edu})}
\and
 Mary Wheeler\thanks{ICES, University of Texas, Austin, TX, USA (\texttt{mfw@ices.utexas.edu})}
}

\maketitle

\begin{abstract}
In this work, we present a novel nonlocal nonlinear coarse grid
approximation using a machine learning algorithm.
We consider unsaturated and two-phase flow problems in heterogeneous and fractured porous media, where mathematical models are formulated as general multicontinuum models. We construct a fine grid approximation using the finite volume method and embedded discrete fracture model. Macroscopic models for these complex nonlinear systems require nonlocal multicontinua approaches, which are developed in earlier works \cite{chung2018nonlinear}. These rigorous techniques require complex local computations, which involve solving local problems in oversampled regions subject to constraints. The solutions of these local problems can be replaced by solving original problem on a coarse (oversampled) region for many input parameters (boundary and source terms) and computing effective properties derived by nonlinear nonlocal multicontinua approaches. The effective properties depend on many variables (oversampled region and the number of continua), thus their calculations require some type of machine learning techniques. In this paper, our contribution is two fold. First, we present macroscopic models and discuss how to effectively compute macroscopic parameters using deep learning algorithms. The proposed method can be regarded as local machine learning and complements our earlier approaches on global machine learning \cite{wang2018deep,wang2019reduced}. We consider a coarse grid approximation using two upscaling techniques with single phase upscaled transmissibilities and nonlocal nonlinear upscaled transmissibilities using a machine learning algorithm. We present results for two model problems in heterogeneous and fractured porous media and show that the presented method is highly accurate and provides fast coarse grid calculations.
\end{abstract}

\section{Introduction}

Mathematical models of the flow and transport problems in heterogeneous and fractured porous media are required to solve large and complex nonlinear systems. Processes in fractured porous media are described by the mixed dimensional coupled system of equations \cite{martin2005modeling, d2012mixed, formaggia2014reduced, lee2001hierarchical, ctene2016algebraic}. Such models can be generalized as a general multicontinuum models similar to the dual porosity/dual permeability approaches \cite{chung2017coupling, vasilyeva2019nonlocal}.

Solving problems in heterogeneous and fractured media requires constructing
grids that resolve all small scale heterogeneity. Numerous model order
reduction techniques have been developed to construct coarse grid
approximations and reduce the computational time of the numerical simulations.
Global model order reduction approaches rely on projection on the important modes space, where the full-order model is carefully study generating POD based basis functions to perform the fast online calculations stage \cite{podporoel2018model, poddeimms2016fast}. Local model reduction techniques are based on constructing local multiscale basis functions to represent the influence of
small scale heterogeneity \cite{ctene2016algebraic, bosma2017multiscale, vasilyeva2019poroel, vasilyeva2018nonlocal}.
One of the widely used  ways is based on the numerical homogenization technique, where effective parameters are calculated in order to construct coarse grid approximations \cite{bakhvalov1984homogenization, sanchez1980non, vasilyeva2018machine}.
The coarse grid parameters are constructed by solving local problems with appropriate boundary conditions. For example, linear boundary conditions or periodicity can be used. The choice of boundary conditions have a strong impact on the accuracy of  results.
In \cite{chen2003coupled}, an interpolated global coarse grid solution is used for performing accurate construction of the upscaled transmissibilities, which involve iterations between global coarse grid model and local fine grid calculations.  In standard upscaling methods, upscaled parameters are obtained independent of any global  problem. However, these approaches lack several features, which are important for rigorous and accurate upscaling. These include the use of multiple continua and oversampled computations.

In the local model order reduction methods,
an oversampling technique and multicontinua concepts are needed to achieve an accuracy independent of physical parameters, such as scales and contrast \cite{chung2018nonfrac,zhao2019analysis}.
For example, the oversampled domain is used for constructing multiscale basis and provide more accurate results in the Multiscale Finite Element Method \cite{hou1997multiscale}. In the Generalized Multiscale Finite Element Method \cite{efendiev2013generalized, chung2014adaptivedg}, a larger domain is used to construct a space of snapshots and solution of the local spectral problem to determine a dominant modes. Note that, the oversampled domain is used for local problem solution and only the interior information is used to define the basis functions.

In recently developed Constrained Energy Minimization and Nonlocal Multicontinuum methods \cite{chung2017cem, chung2018nlmc, vasilyeva2019nlmcmc, vasilyeva2018nlmcperf},  multiscale basis functions are defined in the oversampled domains and constructed via solving local constrained energy minimization problems, where constraints are related to each continuum.
Continuum plays a role of macroscopic parameter.
In \cite{chen2003coupled,durlofsky2007adaptive}, oversampling techniques have been developed in the context of the upscaling procedure, where an interesting local-global upscaling technique is presented for constructing coarse scale approximation for highly heterogeneous porous media. In this method, the coarse grid simulations are iterated with local calculations of the upscaled parameters, where the coarse grid solutions are used to determine the boundary conditions for the local calculation. The local-global upscaling method requires more computation than existing classic upscaling procedures.
Therefore, in the upscaling and multiscale methods oversampled domains are used in two contexts: (1) as extended local domains for more accurate calculations of the coarse grid parameters with fine-scale information about heterogeneous properties and (2) for global or quasi-global information of the solutions that are used, for example, as boundary conditions in local calculations. The second context of the oversampling technique, due to the incorporation of solution information into the local problems leads to the nonlinear equations even in the case of the linear problems.

In this work, we consider flow and transport processes in heterogeneous multicontinuum media and construct coarse and fine grid approximations using a finite volume method with two-point flux approximation.
Due to nonlinear nature of these flows, upscaled parameters are nonlinear functions, which depend on multiple coarse-grid variables defined in oversampled regions. Macroscopic equations use nonlinear nonlocal multicontinuum concept \cite{chung2018nonlinear}. This framework first identifies macroscopic variables in each coarse-grid block and then solves local constraint problems in oversampled regions to compute macroscopic fluxes. The local oversampled computations require solving nonlinear problems with constraints that include the values of macroscopic variables. For example, for two-phase flow simulations, this requires solving two-phase flow problems with known values of pressures and saturations in each coarse-grid block. Each coarse-grid block may contain several macroscale pressures and saturations identified in the first step. Solving local nonlinear constraint problems can be challenging due to large number of nonlinear simulations in oversampled regions. Moreover, computing macroscale fluxes as a function of many variables as a look-up table is nearly impossible. In this work, we propose an efficient algorithm for solving the local problems consisting of original problems with various boundary conditions and using deep learning to train macroscale fluxes. This is a first step in designing computationally efficient and rigorous upscaling methods for nonlinear flows in porous media.

Constructing accurate upscaled transmissibilities for the coarse grid approximation is based on the information about solution (nonlinear transmissibilities).
The presented method is based on the machine learning procedure for fast prediction of the nonlinear transmissibilities, where we construct neural networks that learn dependencies between the coarse grid quantities on the oversampled local domains and upscaled transmissibilities. We use a convolutional neural network and GPU training process to construct a machine learning algorithm \cite{lecun2015deep, krizhevsky2012imagenet}.
For constructing the datasets, we perform local or global calculations of the coarse grid quantities \cite{chen2003coupled}. In the local approach, the upscaled transmissibilities are calculated on the local domain corresponding to the target face, where the fine-scale solution information is used to set boundary conditions. The global approach uses a global fine-scale solution for determining coarse scale parameters.
For training the neural networks, we use a family of problem solutions  for different input conditions.
Note that, we should have many snapshots to capture all input condition variations because accuracy of the machine learning method depends on space of snapshots that is used as a train dataset.
As soon as neural networks trained on the dataset, the fast and accurate calculations can be performed.
%
%
To illustrate method construction and applicability, we considered two model problems: unsaturated flow problem and two-phase filtration problem in heterogeneous and fractured porous media.
The presented method combines accuracy of the fine grid models with fast coarse grid calculations by constructing  machine learning techniques for predicting accurate nonlinear upscaled transmissibilities.

The paper is organized as follows. In Section 2, we present the mathematical model and fine grid approximation. In Section 3, we consider single-phase upscaling for problems in multicontiuum media with variable separation for nonlinear and space dependent variables. Next, we present a novel nonlinear coarse grid approximation using a machine learning algorithm in Section 4. In Section 5, we consider two model problems in two - dimensional formulation and present numerical results, where we consider training of machine learning algorithm and relative errors between the reference fine grid solution and presented method. Finally, we present conclusions.

\section{Model problem with reference fine grid approximation}

As a model problem, we consider two nonlinear problems for fractured and heterogeneous porous media:
\begin{enumerate}
\item Unsaturated flow problem (nonlinear flow problem)
\item Two-phase flow problem (nonlinear transport and flow problem)
\end{enumerate}

We start with the formulation of the mathematical model, where we formulate models for fractured media and generalize it for multicontinuum media. Next, we present a construction of the fine grid approximation using finite volume approximation and embedded fracture model.

\subsection{Unsaturated flow problem}

Mathematical model of the unsaturated flow in porous media described by the Richards' equations \cite{richards1931capillary}
\begin{equation}
\label{uns-m1}
\frac{\partial \Theta}{\partial t}
- \nabla \cdot ( k(x, p) \nabla (p + z)) = f, \quad x \in \Omega,
\end{equation}
where 
$p$ is the pressure head,
$k$ is the unsaturated hydraulic conductivity tensors,
$z$ represent the influence of the gravity to the flow processes,
$\Theta$ is the water content and $f$ refer to source and sink terms.

For fractured porous media, we consider a mixed dimensional formulation of the flow problem \cite{martin2005modeling, d2012mixed, formaggia2014reduced}.
Let $\Omega \in \mathcal{R}^d$ is the $d$ - dimensional domain of the porous matrix, where $d = 2,3$. Fracture network is considered as a $(d-1)$ - dimensional (lower dimensional) domain $\gamma \in \mathcal{R}^{d-1}$ due to small thickness of the fractures compared to the domain sizes.
Then, for unsaturated flow in fractured porous media, we have  the
following coupled system of equations for $p^m$ and $p^f$:
\begin{equation}
\label{uns-m2}
\begin{split}
& \frac{\partial \Theta^m}{\partial t}
- \nabla \cdot ( k^m(x, p^m) \nabla (p^m + z))
+ \sigma^{mf}(x, p^{mf}) (p^m - p^f)= f^m, \quad x \in \Omega, \\
& \frac{\partial \Theta^f}{\partial t}
- \nabla_{\gamma} \cdot ( k^f(x, p^f) \nabla_{\gamma} (p^f + z))
- \sigma^{fm}(x, p^{mf}) (p^m - p^f) = f^f, \quad x \in \gamma,
\end{split}
\end{equation}
where 
$p^m$ and $p^f$ are the pressure head in matrix and fractures;
$k^m$ and $k^f$ are the unsaturated hydraulic conductivity tensors for matrix and fractures;
$z$ represent the influence of the gravity to the flow processes;
$\nabla_{\gamma}$ contains partial derivativies along fracture $\gamma$;
$\Theta^m$ and $\Theta^f$ are the water content for matrix and fracture; and
$f^m$ and $f^f$ refer to source and sink terms.
For transfer term, we have
$\int_{V} \sigma^{mf} (p^m - p^f) dx =
\int_{A} \sigma^{fm} (p^m - p^f) ds$,
$\sigma^{mf} =\sigma/V$,
$\sigma^{fm} = \sigma/A$ ($\sigma = CI \, k^{mf}$ and $k^{mf}(x,p^{mf})$ is the harmonic average between $k^m(x,p^m)$ and $k^f(x,p^f)$) for matrix volume $V$ intersecting with fracture surface and $CI$ is the connectivity index \cite{bosma2017multiscale, ctene2016algebraic}.
As an initial condition, we set $p^{\alpha} = p_0$ ($\alpha = m,f$) and zero flux boundary conditions on $\partial \Omega$ and $\partial \gamma$.

In general, we have following multicontinuum model:
\begin{equation}
\label{uns-m3}
\begin{split}
& \frac{\partial \Theta^{\alpha}}{\partial t}
- \nabla \cdot ( k^{\alpha}(x, p^{\alpha}) \nabla (p^{\alpha} + z))
+ \sum_{\beta} \sigma^{\alpha \beta}(x,p^{\alpha \beta}) (p^{\alpha} - p^{\beta}) = f^{\alpha},
\end{split}
\end{equation}
where $\alpha = 1,...,M$ and $M$ is the number of continuum.

Let $c^{\alpha}(x, p^{\alpha}) = \partial \Theta^{\alpha} / \partial p$ 
therefore we have following coupled system of nonlinear parabolic equations
\begin{equation}
\label{uns-m4}
c^{\alpha}(x, p^{\alpha}) \frac{\partial p^{\alpha}}{\partial t}
 - \nabla \cdot (k^{\alpha}(x, p^{\alpha}) \nabla p^{\alpha})
 + \sum_{\beta}
 \sigma^{\alpha \beta}(x,p^{\alpha \beta}) (p^{\alpha} - p^{\beta})
 = q^{\alpha},
\end{equation}
where $q^{\alpha} = f^{\alpha} + \nabla \cdot (k^{\alpha}(x, p^{\alpha}) z)$, $c^{\alpha}(p^{\alpha})$ and $k^{\alpha}(x, p^{\alpha})$ are the nonlinear coefficient ($\alpha = 1,...,M$).

For the approximation on the fine grid, we use structured grids with embedded discrete fracture model (EDFM) \cite{lee2001hierarchical, ctene2016algebraic}.
Let $\mathcal{T}^h$ denote a structured fine grid of the porous matrix domain  $\Omega$ and $\mathcal{G}^h$ denote a fine grid of the fracture domain  $\gamma$
\[
\mathcal{T}^h = \cup_{i=1}^{N^{m,h}} \varsigma_i, \quad
\mathcal{G}^h = \cup_{l=1}^{N^{f,h}} \iota_l,
\]
where $\varsigma_i$ and $\iota_l$ are the cell of the matrix and fractures fine grids,
$N^{m,h}$ is the number of cells in $\mathcal{T}^h$,
$N^{f,h}$ is the number of cell related to fracture mesh $\mathcal{G}^h$.
Therefore, for finite volume approximation we have
\[
\begin{split}
& c^m_i \frac{ p_i^m - \check{p}_i^m }{\tau} |\varsigma_i|
 + \sum_{j}  u_{ij}^{mm} +  \sum_{l} u^{mf}_{il}
 =  q^m_i   |\varsigma_i|,
 \quad \forall i = 1,..., N^{m,h} \\
& c^f_l \frac{ p_l^f - \check{p}_l^f}{\tau}  |\iota_l|
+ \sum_{n}  u_{ln}^{ff}
+ \sum_{i} u^{fm}_{il}
 =  q^f_l  |\iota_l|,
 \quad \forall l = 1,.., N^{f, h},
\end{split}
\]
where
$p_i^m$ and $p_l^f$ are pressure of matrix and fracture continuum in cells $\varsigma_i$ and $\iota_l$,
$|\varsigma_i|$ and $|\iota_l|$ are the volume of cells.
Here, we use an implicit scheme for time discretization, where
$\check{p}^m_i$ and  $\check{p}^f_l$ are the solutions from previous time step and $\tau$ is the given time step \cite{vasilyeva2019poroel, vasilyeva2018nonlocal}.

For approximation of the flux in matrix ($u^{mm}$) and fracture ($u^{ff}$) continuum
\[
u^{mm} =  - k^m(x, p^m) \nabla p^m, \quad
u^{ff} = - k^f(x, p^f) \nabla_{\gamma} p^f,
\]
we use a classic two point flux approximation (TPFA)
\[
u_{ij}^{mm} = u^{mm} \cdot n |_{E_{ij}}  =
T^{mm}_{ij}(p^m_i, p^m_j) (p^m_i - p^m_j), \quad
u_{ln}^{ff} =  u^{ff} \cdot n |_{e_{ln}}  =
 T^{ff}_{ij}(p^f_i, p^f_l) (p^f_i - p^f_j),
\]
where
$T_{ij}^{mm} = k^m_{ij} |E_{ij}|/d_{ij}$,
$T_{ln}^{ff} = k^f_{ln}/\Delta_{ln}$,
$k^m_{ij} = (k_i^m(p_i^m) + k_j^m(p_j^m))/2$,
$k^f_{ln} = (k_l^f(p_l^f) + k_n^f(p_n^f))/2$,
$E_{ij}$ and $e_{ln}$ are the interface between two cells,
$|E_{ij}|$ is the length of face between cells $\varsigma_i$ and $\varsigma_j$, $d_{ij}$ is the distance between midpoint of cells $\varsigma_i$ and $\varsigma_j$,
$\Delta_{ln}$ is the distance  between midpoint of cells $\iota_l$ and $\iota_n$.

For the flux  between matrix and fracture continuum
\[
u^{mf} = \sigma^{mf}(x, p^{mf}) (p^m - p^f), \quad
u^{fm} = \sigma^{fm}(x, p^{mf}) (p^f - p^m),
\]
we follow EDFM and have the following approximation
\[
u^{mf}_{il} = -u^{fm}_{il} = T^{mf}_{il} (p^m_i - p^f_l),
\]
where $T^{mf}_{il} = \sigma_{il}$ with  $\sigma_{il} =  CI_{il} k^{mf}_{il}$ if $\iota_l \subset \varsigma_i$ and zero else ($CI_{il}$ is the connectivity index from \cite{lee2001hierarchical, ctene2016algebraic, hkj12} that proportional to the distance and area of the intersection between the fracture cell $\iota_l$ and porous matrix cell $\varsigma_i$).

Therefore, we have following fine grid approximation
\begin{equation}
\label{uns-f1}
\begin{split}
& c^m_i \frac{ p_i^m - \check{p}_i^m }{\tau} |\varsigma_i|
 + \sum_{j}  T_{ij}^{mm} (p^m_i, p^m_j)  (p_i^m - p_j^m)
 + \sum_l T^{mf}_{il}(p^m_i, p^f_l) (p_i^m - p_l^f )
 =  q^m_i   |\varsigma_i|,
 \quad \forall i = 1,..., N^{m,h} \\
& c^f_l \frac{ p_l^f - \check{p}_l^f}{\tau}  |\iota_l|
+ \sum_{n}  T_{ln}^{ff}(p^f_l, p^f_n)  (p_l^f - p_n^f)
- \sum_i T^{mf}_{il}(p^m_i, p^f_l) (p_i^m - p_l^f )
 =  q^f_l  |\iota_l|,
 \quad \forall l = 1,.., N^{f, h},
\end{split}
\end{equation}
or for  general multicontinuum case \cite{chung2017coupling, vasilyeva2019nonlocal}, we have
\begin{equation}
\label{uns-f4}
 c^{\alpha}_i \frac{ p^{\alpha}_i - \check{p}^{\alpha}_i }{\tau} |\varsigma_i|
 + \sum_j
 T^{\alpha \alpha}_{ij} (p^{\alpha}_i, p^{\alpha}_j)  (p^{\alpha}_i - p^{\alpha}_j)
 + \sum_{\beta} \sum_l
 T^{\alpha \beta}_{il}(p^{\alpha}_i, p^{\beta}_l)(p^{\alpha}_i - p^{\beta}_l )
 =  q^{\alpha}_i   |\varsigma_i|, \quad \forall i = {1,..., N^{\alpha, h} },
\end{equation}
where $\alpha = 1,...,M$ and  $M$ is the number of continua.

\subsection{Two-phase flow problem }

Mathematical model of the two-phase flow problem in porous media contains a conservation law and Darcy's law \cite{helmig1997multiphase}. For the case with incompressible fluid and rock and without gravitational and capillary forces, we have
\begin{equation}
\label{tp-m1}
\begin{split}
\phi \frac{ \partial s}{\partial t}
- \nabla \cdot ( \lambda^w(s) k(x) \nabla p) = q^w, \quad x \in \Omega, \\
- \nabla \cdot (\lambda(s) k(x) \nabla p) = q, \quad x \in \Omega,
\end{split}
\end{equation}
where $s = s^w$ is the saturation of the wetting phase, $p$ is the pressure,
$q = q^w + q^n $, $q^w$ and $q^n$ are the source/sink of wetting and nonwetting phases,
$\lambda^i = {k_{ri}(s)}/{\mu_i}$, $\lambda = \lambda^n + \lambda^w$,
$\phi$, $k$ are the porosity and permeability, $\mu_i$ and $k_{ri}$ are the viscosity and relative permeability for $i$-phase ($i = n, w$).

For the fractured porous media, we consider the mixed dimensional mathematical model for two-phase flow problem
\begin{equation}
\label{tp-m2}
\begin{split}
\phi^m \frac{ \partial s^m}{\partial t}
- \nabla \cdot ( \lambda^w(s^m) k^m(x) \nabla p^m)
+ \lambda^w(s^{mf}) \sigma^{mf}(x) (p^m - p^f)
= q^{w,m}, \quad x \in \Omega, \\
- \nabla \cdot (\lambda(s^m)  k^m(x) \nabla p^m)
+ \lambda(s^{mf}) \sigma^{mf}(x) (p^m - p^f)
= q^m, \quad x \in \Omega, \\
\phi^f \frac{ \partial s^f}{\partial t}
- \nabla_{\gamma} \cdot ( \lambda^w(s^f) k^f(x) \nabla_{\gamma} p^f)
-  \lambda^w(s^{mf}) \sigma^{fm}(x) (p^m - p^f)
= q^{w,f}, \quad x \in \gamma, \\
- \nabla_{\gamma} \cdot (\lambda(s^f)  k^f(x) \nabla_{\gamma} p^f)
+ \lambda(s^{mf}) \sigma^{fm}(x) (p^m - p^f)
= q^f, \quad x \in \gamma, \\
\end{split}
\end{equation}
where $s^m$ and $s^f$ are the saturation in porous matrix and in fractures;
$p^m$ and $p^f$ are the pressure in matrix and in fractures;
$\phi^m$ and $\phi^f$ are the porosity for matrix and fracture continuum;
$k^m$ and $k^f$ are the absolute permeability of matrix and fractures;
and $q^{\alpha} = q^{w,\alpha} + q^{n,\alpha} $, $q^{w,\alpha}$ and $q^{n,\alpha}$ are the source/sink of wetting and nonwetting phases for continua $\alpha = m,f$.
On the fine grid, we suppose $\lambda^{i,m} = \lambda^{i,f} = \lambda^i$ for $i = n, w$ and $\lambda = \lambda^n + \lambda^w$, where
$\lambda^i = {k_{ri}(s^{\alpha})}/{\mu_i}$, $\mu_i$ is the viscosity and $k_{ri}$ is relative permeability for $i$-phase that depends on flux direction.  In general, realtive permeability functions can be different for each continua and for flux between them.
For transfer term, similarly to the previous model for unsaturated flow, we have
$\sigma^{mf} =\sigma/V$,
$\sigma^{fm} = \sigma/A$ ($\sigma = CI \, k^{mf}$ and $k^{mf}$ is the harmonic average between $k^m$ and $k^f$) for matrix volume $V$ intersecting with fracture surface and $CI$ is the connectivity index \cite{bosma2017multiscale, ctene2016algebraic}.
As an initial condition, we set $s^{\alpha} = s_0(x)$  ($\alpha = m,f$). For boundary conditions, we set zero flux on $\partial \Omega$ and $\partial \gamma$.

For the general multicontinuum case, we can write
\begin{equation}
\label{tp-m3}
\begin{split}
\phi^{\alpha} \frac{ \partial s^{\alpha}}{\partial t}
- \nabla \cdot ( \lambda^w(s^{\alpha}) k^{\alpha}(x) \nabla p^{\alpha})
+ \sum_{\beta}
\lambda^w(s^{\alpha \beta}) \sigma^{\alpha \beta}(x) (p^{\alpha} - p^{\beta})
= q^{w,\alpha}, \\
- \nabla \cdot (\lambda (s^{\alpha})  k^{\alpha}(x) \nabla p^{\alpha})
+ \sum_{\beta}
\lambda(s^{\alpha \beta}) \sigma^{\alpha \beta}(x) (p^{\alpha} - p^{\beta})
= q^{\alpha},
\end{split}
\end{equation}
where $\alpha = 1,...,M$ and $M$ is the number of continuum.

Similarly to the previous model problem of unsaturated flow, we use structured grids and construct a finite volume approximation with embedded discrete fracture model (EDFM) for approximation of the two-phase flow problem.
We use same fine grid for the porous matrix domain ($\mathcal{T}^h =  \cup_{i=1}^{N^{m,h}} \varsigma_i$) and the fracture domain ($\mathcal{G}^h = \cup_{l=1}^{N^{f,h}} \iota_l$) with cells  $\varsigma_i$ and $\iota_l$,
$N^{m,h}$ is the number of cells in $\mathcal{T}^h$,
$N^{f,h}$ is the number of cell related to fracture mesh $\mathcal{G}^h$.
For approximation by time, we use IMPES (implicit pressure explicit saturation) scheme and a finite volume approximation by space
\[
\begin{split}
\phi^m_i \frac{s^m_i - \check{s}^m_i}{\tau} |\varsigma_i| + &
\sum_j  u_{ij}^{w,mm} + \sum_l u^{w, mf}_{il}
=  (1 - f^{w,m}_i) q_i^{w,m} |\varsigma_i|,  \\
&
\sum_j  u_{ij}^{mm}  +  \sum_l u^{mf}_{il} = q_i^m |\varsigma_i|,
\\
 \phi^f_l \frac{ s_l^f - \check{s}_l^f}{\tau}  |\iota_l| + &
\sum_n  u_{ln}^{w, ff} + \sum_i u^{w, fm}_{il}
 =   (1 - f^{w,f}_l)  q_l^{w,f}  |\iota_l|, \\
& \sum_n  u_{ln}^{ff} + \sum_i u^{fm}_{il} =  q_l^f  |\iota_l|,
\end{split}
\]
where
$s_i^m$ and $s_l^f$ are saturation of matrix and fracture continuum in cells $\varsigma_i$ and $\iota_l$,
$f^{w,\alpha}_i = {\lambda^w(s^{\alpha}_i)}/{\lambda(s^{\alpha}_i)}$,
$|\varsigma_i|$ and $|\iota_l|$ are the volume of cells.
Here, we use an implicit scheme for time discretization, where
$\check{s}^m_i$ and  $\check{s}^f_l$ are the solutions from previous time step and $\tau$ is the given time step.

For the fluxes
\[
u^{mm} =  - \lambda(s^m) k^m(x) \nabla p^m, \quad
u^{w, mm} =  - \lambda^w(s^m) k^m(x) \nabla p^m,
\]\[
u^{ff} = - \lambda(s^f) k^f(x) \nabla_{\gamma} p^f,\quad
u^{w, ff} = - \lambda^w(s^f) k^f(x) \nabla_{\gamma} p^f,
\]\[
u^{mf} = \lambda(s^{mf}) \sigma^{mf}(x) (p^m - p^f), \quad
u^{w, mf} = \lambda^w(s^{mf}) \sigma^{mf}(x) (p^m - p^f),
\]\[
u^{fm} = \lambda(s^{mf}) \sigma^{fm}(x)  (p^f - p^m), \quad
u^{w, fm} = \lambda^w(s^{mf}) \sigma^{fm}(x)  (p^f - p^m),
\]
we have following approximations
\[
u_{ij}^{mm} = u^{mm} \cdot n |_{E_{ij}}  =
T^{mm}_{ij}  (p^m_i - p^m_j), \quad
u_{ij}^{w,mm} = u^{w,mm} \cdot n |_{E_{ij}}  =
T^{w,mm}_{ij}  (p^m_i - p^m_j),
\]\[
u_{ln}^{ff} =  u^{ff} \cdot n |_{e_{ln}}  =
 T^{ff}_{ij}(p^f_i - p^f_j), \quad
u_{ln}^{w,ff} =  u^{w,ff} \cdot n |_{e_{ln}}  =
 T^{w,ff}_{ij} (p^f_i - p^f_j),
\]\[
u^{mf}_{il} = -u^{fm}_{il} = T^{mf}_{il} (p^m_i - p^f_l), \quad
u^{w, mf}_{il} = -u^{w, fm}_{il} = T^{w, mf}_{il} (p^m_i - p^f_l),
\]
where
\[
T_{ij}^{\alpha \beta}
= T_{ij}^{\alpha \beta} (s^{\alpha}_i, s^{\beta}_j)
= \lambda(s^{\alpha \beta}_{ij}) W_{ij}^{\alpha \beta}, \quad
T_{ij}^{w, \alpha \beta}
= T_{ij}^{w, \alpha \beta} (s^{\alpha}_i, s^{\beta}_j)
= \lambda^w (s^{\alpha \beta}_{ij}) W_{ij}^{\alpha \beta},  \quad \alpha, \beta = m,f
\]
and
$W_{ij}^{mm} = k^m_{ij} |E_{ij}|/d_{ij}$,
$W_{ln}^{ff} = k^f_{ln}/\Delta_{ln}$,
$W_{il}^{mf} = \sigma_{il}$  with  $\sigma_{il} =  CI_{il} k^m_{il}$ if $\iota_l \subset \varsigma_i$ and zero else.  Here $|E_{ij}|$ is the length of face between cells $\varsigma_i$ and $\varsigma_j$, $d_{ij}$ is the distance between midpoint of cells $\varsigma_i$ and $\varsigma_j$, $\Delta_{ln}$ is the distance between points $l$ and $n$, $CI_{il}$ is the connectivity index from \cite{lee2001hierarchical, ctene2016algebraic, hkj12} that proportional to the distance and area of the intersection between the fracture cell $\iota_l$ and porous matrix cell $\varsigma_i$.

Therefore, we have following discrete problem on the fine grid
\begin{equation}
\label{tp-m4}
\begin{split}
\phi^m_i \frac{s^m_i - \check{s}^m_i}{\tau} |\varsigma_i|
+ & \sum_{j}  T_{ij}^{w,mm} (\check{s}^m_i, \check{s}^m_j)  (p_i^m - p_j^m)
+  \sum_l T^{w, mf}_{il}(\check{s}^m_i, \check{s}^f_l) (p_i^m - p_l^f )
=   (1 - f^{w,m}_i) q_i^{w,m} |\varsigma_i|,  \\
& \sum_{j}  T_{ij}^{mm} (\check{s}^m_i, \check{s}^m_j)  (p_i^m - p_j^m)
+  \sum_l T^{mf}_{il}(\check{s}^m_i, \check{s}^f_l) (p_i^m - p_l^f )
= q_i^m |\varsigma_i|,
\\
 \phi^f_l \frac{ s_l^f - \check{s}_l^f}{\tau}  |\iota_l|
+ & \sum_{n}  T_{ln}^{w, ff}(\check{s}^f_l, \check{s}^f_n)  (p_l^f - p_n^f)
- \sum_i T^{w, mf}_{il}(\check{s}^m_i, \check{s}^f_l) (p_i^m - p_l^f )
 =  (1 - f^{w,f}_l) q_l^{w,f}  |\iota_l|, \\
& \sum_{n}  T_{ln}^{ff}(\check{s}^f_l, \check{s}^f_n)  (p_l^f - p_n^f)
- \sum_i T^{mf}_{il}(\check{s}^m_i, \check{s}^f_l) (p_i^m - p_l^f )
 =  q_l^f  |\iota_l|,
\end{split}
\end{equation}
where $i = 1,..., N^{m,h}$ and $ l = 1,.., N^{f, h}$.

For approximation of the $\lambda^w(s^{\alpha \beta}_{ij})$, we use an upwind scheme
\[
\lambda^w(s^{\alpha \beta}_{ij})  =
\begin{cases}
\lambda^w(s^{\alpha}_i), 	& \mbox{if } T_{ij}^{\alpha \beta} (\check{s}^{\alpha}_i, \check{s}^{\beta}_j)  (p_i^{\alpha} - p_j^{\beta}) > 0 \\
\lambda^w(s^{\beta}_j), 	& else,
\end{cases}.
\]
and $ \lambda(s^{\alpha \beta}_{ij})$ is the harmonic average between $\lambda(s^{\alpha}_i)$ and $\lambda(s^{\beta}_j)$.

For the general multicontinuum model, we have
\begin{equation}
\label{tp-m5}
\begin{split}
\phi^{\alpha}_i \frac{s^{\alpha}_i - \check{s}^{\alpha}_i}{\tau} |\varsigma_i|
+ & \sum_{j}  T_{ij}^{w,\alpha \alpha} (\check{s}^{\alpha}_i, \check{s}^{\alpha}_j)  (p_i^{\alpha} - p_j^{\alpha})
+  \sum_{\beta} \sum_l
T^{w, \alpha \beta}_{il}(\check{s}^{\alpha}_i, \check{s}^{\beta}_l) (p_i^{\alpha} - p_l^{\beta} )
=  (1 - f^{w,\alpha}_i)  q_i^{w,\alpha} |\varsigma_i|,  \\
& \sum_{j}  T_{ij}^{\alpha \alpha} (\check{s}^{\alpha}_i, \check{s}^{\alpha}_j)  (p_i^{\alpha} - p_j^{\alpha})
+  \sum_{\beta} \sum_l
T^{\alpha \beta}_{il}(\check{s}^{\alpha}_i, \check{s}^{\beta}_l) (p_i^{\alpha} - p_l^{\beta} )
= q_i^{\alpha} |\varsigma_i|, \\
\end{split}
\end{equation}
where $\alpha = 1,...,M$.

\section{Coarse grid upscaled model}

Solution of the problems in heterogeneous and fractured media require
the construction of the grids that resolve all small scale heterogeneity.
One of the widely used  way is based on the numerical homogenization
technique, where effective parameters are calculated in order
to construct a coarse grid approximations.
The coarse grid parameters are constructed by a solution of
the local problems with appropriate boundary conditions.

\begin{figure}[h!]
\centering
\includegraphics[width=0.49\linewidth]{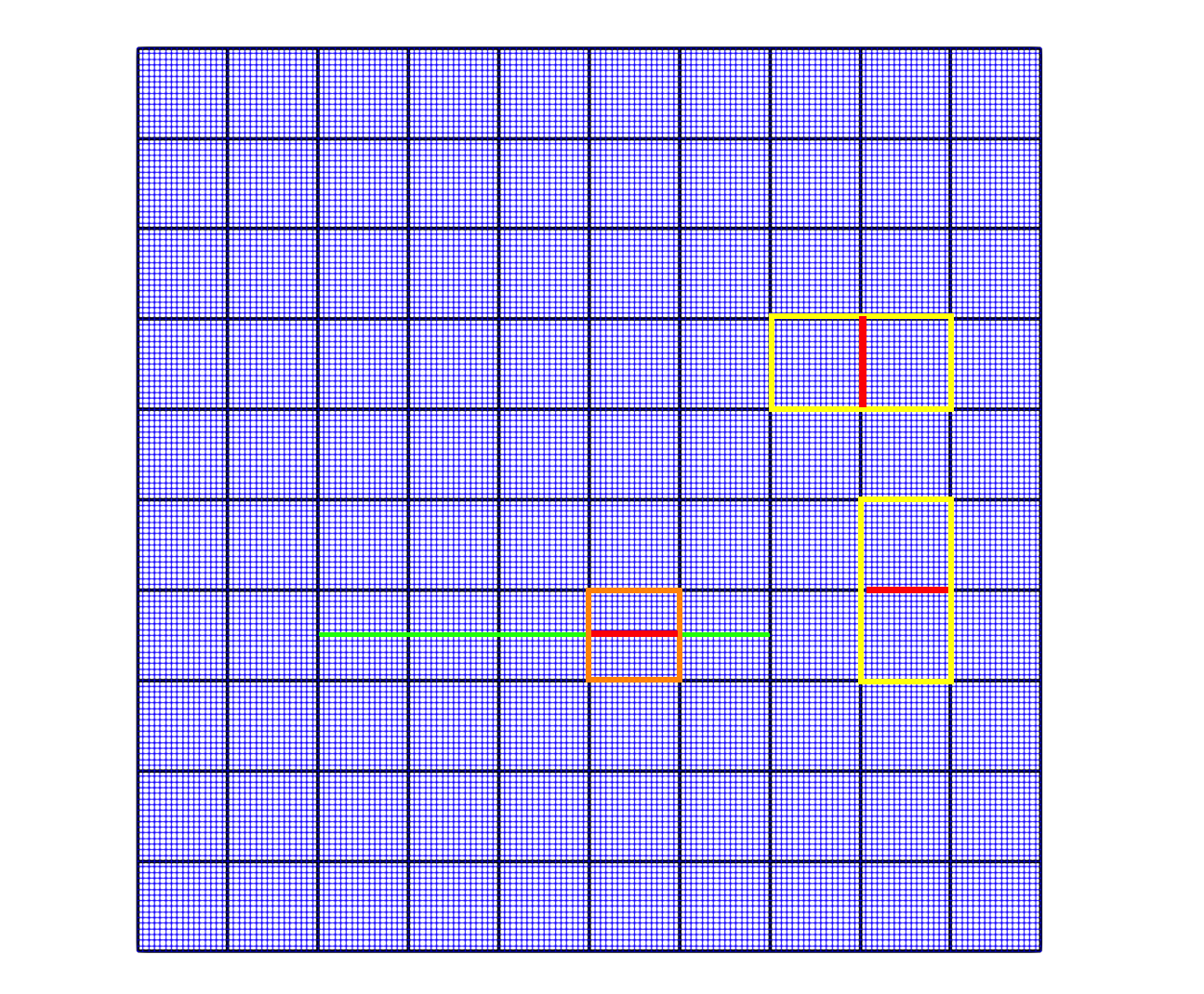}
\includegraphics[width=0.49\linewidth]{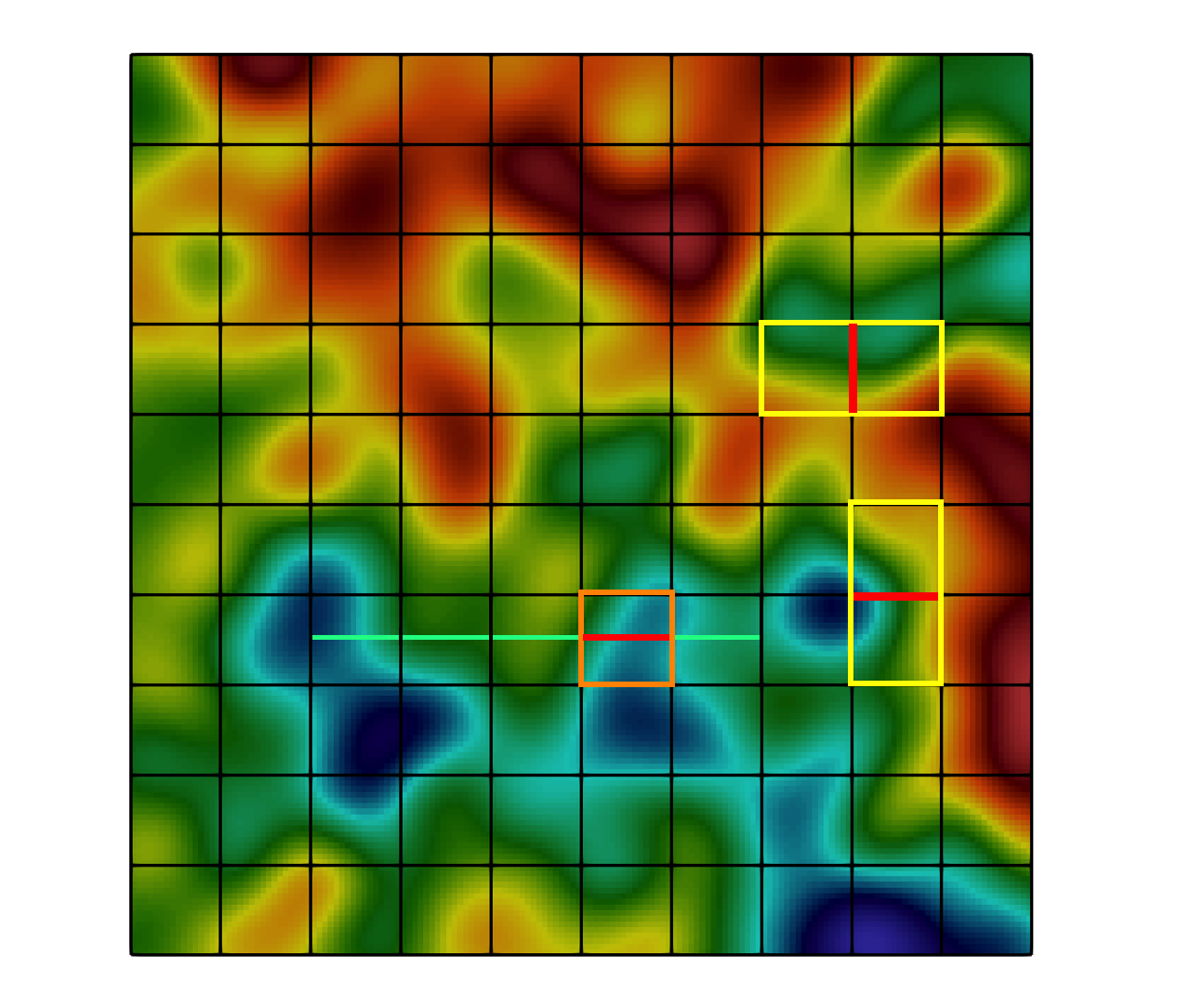}
\caption{Coarse grid, heterogeneous properties and local domains.
Left: coarse grid with fine grid and local domains.
Right: heterogeneous permeability with coarse cells and local domains.
Coarse grid  $\mathcal{T}^H$ (black color), fine grid $\mathcal{T}^h$ (blue color), fracture $\gamma$ (green), local domain $\omega_{ij}$ (yellow), coarse edge $E_{ij}$ (red) }
\label{fig:meshc}
\end{figure}

Let $\mathcal{T}^H$ be a structured coarse mesh of the computational domain $\Omega$
\[
\mathcal{T}^H = \cup_{i=1}^{N^H} K_i,
\]
where $N^H$ is the number of the coarse grid cells, $K_i$ is the quadrilateral coarse cell and $i$ is the coarse grid cell index \cite{vasilyeva2018machine, vasilyeva2019convolutional}.
Form of the coarse grid upscaled model is similar to the fine grid model with finite volume approximation, where coarse grid transmissibilities $T_{ij}^{UP}$ are calculated by a solution of the local problems that take into account fine grid resolution of the heterogeneous permeability  (see Figure \ref{fig:meshc}).

We let $E_{ij}$ be the coarse grid face and we define the neighborhood (local domain) by
\[
\omega_{ij} = K_i \cup K_j, \quad K_i, K_j \in \mathcal{T}^H,
\]
where $\omega_{ij}$ is a union of two coarse cells, when $E_{ij}$ lies in the interior of the domain $\Omega$. For the edges on the boundary, we will use a no flux boundary conditions and therefore not need to calculate of the upscaled transmissibilities.
For calculation of the upscaled transmissibilities $T_{ij}^{UP}$ for coarse face $E_{ij}$, we solve local problems for nonperiodic heterogeneous fractured media with linear boundary conditions in $\omega_{ij}$. 
For the fractured/multicontinuum media, we use a similar approach for calculation of the coarse grid transmissibilities $T_{ij}^{\alpha \beta, UP}$. Details of the calculations, we present below for each problem.


\subsection{Nonlinear flow problem}

We start with single-phase upscaling and suppose that
\[
k^{\alpha}(x, p^{\alpha})
= k_r(p^{\alpha}) k_s^{\alpha}(x),  \quad
\sigma^{\alpha \beta}(x, p^{\alpha \beta})
= \sigma_r(p^{\alpha \beta}) \sigma^{\alpha \beta}_s(x), \quad
\quad \alpha = m, f,
\]
where, in general, $k_r$ and $\sigma_r$ can be different for each continuum $\alpha$, but in this work, we assume similar relationships, for simplicity.

Let $\mathcal{T}^H$ denote a structured coarse grid of the porous matrix domain  $\Omega$ and $\mathcal{G}^H$ denote a coarse grid of the fracture domain  $\gamma$
\[
\mathcal{T}^H = \cup_{i=1}^{N^{m,H}} K_i, \quad
\mathcal{G}^H = \cup_{l=1}^{N^{f,H}} \gamma_l,
\]
where $K_i$ and $\gamma_l$ are the cell of the matrix and fractures fine grids,
$N^{m,H}$ is the number of coarse cells in $\mathcal{T}^H$,
$N^{f,H}$ is the number of coarse cells related to $\mathcal{G}^H$.
On the coarse grid for equation \eqref{uns-m4}, we have the following
discrete problem for $\overline{p} = (\overline{p}^m, \overline{p}^f)$
\begin{equation}
\label{uns-up1}
\begin{split}
\overline{c}^m_i
\frac{ \overline{p}_i^m - \check{\overline{p}}_i^m }{\tau} |K_i|
 + \sum_{j}  {T}_{ij}^{mm, UP}  (\overline{p}_i^m - \overline{p}_j^m)
 + \sum_l {T}^{mf, UP}_{il} (\overline{p}_i^m - \overline{p}_l^f )
 =  \overline{q}^m_i   |K_i|,  \\
\overline{c}^f_l \frac{ \overline{p}_l^f - \check{\overline{p}}_l^f}{\tau}  |\gamma_l|
+ \sum_{n}  {T}_{ln}^{ff, UP} (\overline{p}_l^f - \overline{p}_n^f)
- \sum_i {T}^{mf, UP}_{il} (\overline{p}_i^m - \overline{p}_l^f )
 =  \overline{q}^f_l  |\gamma_l|,
\end{split}
\end{equation}
where $l = 1,.., N^{f, H}$ and $i = 1,..., N^{m,H}$.

In general for multicontinuum model, we have
 \begin{equation}
\label{uns-up2}
\begin{split}
\overline{c}^{\alpha}_i
\frac{ \overline{p}_i^{\alpha} - \check{\overline{p}}_i^{\alpha} }{\tau} |K_i|
 + \sum_{j}
 {T}_{ij}^{\alpha \alpha, UP}  (\overline{p}_i^{\alpha} - \overline{p}_j^{\alpha})
 + \sum_{\beta}   \sum_l
 {T}^{\alpha \beta, UP}_{il} (\overline{p}_i^{\alpha} - \overline{p}_l^{\beta} )
 =  \overline{q}^{\alpha}_i   |K_i|,
\end{split}
\end{equation}
where $\overline{c}^{\alpha}_i \approx  \frac{1}{|K_i|} \int_{K_i} c^{\alpha}(\overline{p}_i^{\alpha}) \, dx $  and
\begin{equation}
\label{uns-up3}
{T}^{\alpha \beta, UP}_{ij}(\overline{p}_i^{\alpha}, \overline{p}_j^{\beta})
=
k_r(\overline{p}_{ij}^{\alpha \beta})
W^{\alpha \beta, UP}_{ij}, \quad \alpha = m, f
\end{equation}
and $W^{\alpha \beta, UP}_{ij}$ is the precalculated effective transmissibilities.

For the calculation of the upscaled transmissibilities for the porous matrix,
we solve the following local problems in each $\omega_{ij}$  (see Figure \ref{fig:meshc}, where local domain is depicted by a yellow color)
\begin{equation}
\label{uns-loc-Tup}
- \nabla \left(k_s^{m}(x) \nabla \psi^l \right) = 0, \quad x \in \omega_{ij},
\end{equation}
with boundary conditions
\[
\begin{split}
&\psi^l = 1, \quad x \in \Gamma_{ij}^1,\quad
\psi^l = 0, \quad x \in \Gamma_{ij}^2,\\
-& k_s^{m}(x) \frac{\partial \psi^l}{\partial n} = 0,
\quad x \in \partial \omega_{ij}/ (\Gamma_{ij}^1 \cup \Gamma_{ij}^2).\\
\end{split}
\]
In this work, we consider two-dimensional problems with $x = (x_1, x_2)$. Therefore, we solve two local problems for $\psi^l$, $l = 1, 2$.
For $\psi^1$, boundaries $\Gamma_{ij}^1$ and $\Gamma_{ij}^2$ are the left and right boundaries of the domain $\omega_{ij}$, respectively.
For $\psi^2$, boundaries $\Gamma_{ij}^1$ and $\Gamma_{ij}^2$ are the top and bottom boundaries of the domain $\omega_{ij}$, respectively.
Note that, another boundary conditions can be applied for local problems.

Therefore, for calculations  $W_{ij}^{mm,UP}$ in \eqref{uns-up3}, we solve following discrete problem for finite volume approximation up to fine grid resolution
\[
\sum_j W_{ij} (\psi^l_i - \psi^l_j) = 0,
\]
with appropriate boundary conditions.

After solution of the local problems in $\omega_{ij}$, we calculate
upscaled transmissibility for the porous matrix (see Figure \ref{fig:meshc}, where interface $E_{ij}$ is depicted by a red color)
\begin{equation}
\label{uns-Tup}
W^{mm, UP}_{ij} =
 \frac{
 \sum_{r,n} W_{rn} (\psi^l_r - \psi^l_n) }{\overline{\psi}^l_i - \overline{\psi}^l_j},
\end{equation}
where $r,n$ are the fine cells around coarse face $E_{ij}$, 
$\overline{\psi}^l_i$ and $\overline{\psi}^l_j$ are the mean values in coarse cells $K_i$ and $K_j$. We use $\psi^l$ with $l=1$ for all vertical edges and $l=2$ for horizontal edges.
For the fracture continuum, we suppose that $k^f = \const$ and therefore set $W_{ln}^{ff, UP} = k^f/d_{ln}$ ($d_{ln}$ is the distance between points $l$ and $n$).  

In this work, we suppose $k^f = const$, and therefore for the calculations of the coarse grid transmissibility between coarse grid fracture cells $\gamma_l$ and $\gamma_n$, we have
$T^{ff, UP}_{ln} = k^f/\Delta_{ln}$,  where $\Delta_{ln}$ is the distance  between midpoint of cells $\gamma_l$ and $\gamma_n$.

Let $\omega^{mf}_{il} = \{K_i : K_i \cup \gamma_l \neq \emptyset, \,  \gamma_l \in \mathcal{G}^H, \, K_i \in \mathcal{T}^H\}$ be the local domain for calculation of the $W^{mf}_{il}$ in \eqref{uns-up3} (see Figure \ref{fig:meshc}, where local domain is depicted by a orange color).
For the calculation of the upscaled transmissibility between porous matrix and fracture, we solve local problems in  $\omega^{mf}_{il}$
\begin{equation}
\label{uns-loc-Tup-fm}
c^m \frac{\partial \phi}{\partial t}
- \nabla \left(k_s^{m}(x) \nabla \phi \right) + \sigma^{mf}_s(x)(\phi - \phi^f) = 0, \quad x \in \omega^{mf}_{il},
\end{equation}
where $\phi^f = 1$ on  $\gamma_l $ with zero flux boundary conditions on $\partial  \omega^{mf}_{il}$.
Therefore, we solve following discrete system for finite volume approximation up
to fine grid resolution
\[
c^m_i \frac{\phi_i - \check{\phi}_i}{\tau} |\varsigma_i|
+ \sum_j W_{ij} (\phi_i - \phi_j)
+ \sum_l W^{mf}_{il} (\phi_i - \phi^f_l) = 0,
\]
until $|\phi_i - \check{\phi}_i| > \epsilon$ and find upscaled matrix-fracture transmissibility using final time step solution
\begin{equation}
\label{uns-Tup-fm}
W^{mf}_{il} =
\frac{
\sum_{r,n} W^{mf}_{rn} (\phi_r - \phi^f_n) }{\overline{\phi}_i - \overline{\phi}^f_l},
\end{equation}
where $r$ are the cell that contains fracture, 
$\overline{\phi}_i$ and $\overline{\phi}^f_l$ are the mean values in coarse cells $K_i$ and in fracture $\gamma_l$ (see Figure \ref{fig:meshc}, where interface $\gamma_l$ is depicted by a red color).

Note that, there exist different approaches for calculation of the effective transmissibilities, for example, based on the different boundary conditions for local problems, using oversampled domains and the construction of look-up table for interpolation of the nonlinear dependence. In this work, for calculating the  upscaled transmissibilities, we use the simplest classic approach. The main goal of the paper is the construction of the novel highly accurate  nonlinear upscaled coarse grid approximations using machine learning techniques.

\subsection{Nonlinear flow and transport problem }

On the coarse grid for equation \eqref{tp-m2}, we have following discrete problem for $\overline{p}^m, \overline{p}^f, \overline{s}^m$ and  $\overline{s}^f$
\begin{equation}
\label{tp-up1}
\begin{split}
\overline{\phi}^m_i
\frac{\overline{s}^m_i - \check{\overline{s}}^m_i}{\tau} |K_i|
+ & \sum_{j}
T^{w,mm, UP}_{ij} (\check{\overline{s}}^m_i, \check{\overline{s}}^m_j)  (\overline{p}_i^m - \overline{p}_j^m)
+  \sum_l
T^{w, mf, UP}_{il}(\check{\overline{s}}^m_i, \check{\overline{s}}^f_l) (\overline{p}_i^m - \overline{p}_l^f )
=   (1 - \overline{f}^{w,m}_i) \overline{q}_i^{w,m} |K_i|,  \\
& \sum_{j}
T^{mm, UP}_{ij} (\check{\overline{s}}^m_i, \check{\overline{s}}^m_j)  (\overline{p}_i^m - \overline{p}_j^m)
+  \sum_l
T^{mf, UP}_{il}(\check{\overline{s}}^m_i, \check{\overline{s}}^f_l) (\overline{p}_i^m - \overline{p}_l^f )
= \overline{q}_i^m |K_i|,
\\
\overline{\phi}^f_l \frac{ \overline{s}_l^f - \check{\overline{s}}_l^f}{\tau}  |\gamma_l|
+ & \sum_{n}
T^{w, ff, UP}_{ln} (\check{\overline{s}}^f_l, \check{\overline{s}}^f_n)  (\overline{p}_l^f - \overline{p}_n^f)
- \sum_i
T^{w, mf, UP}_{il}(\check{\overline{s}}^m_i, \check{\overline{s}}^f_l) (\overline{p}_i^m - \overline{p}_l^f )
 =  (1 - \overline{f}^{w,f}_i) \overline{q}_l^{w,f}  |\gamma_l|, \\
& \sum_{n}
T^{ff, UP}_{ln} (\check{\overline{s}}^f_l, \check{\overline{s}}^f_n)  (\overline{p}_l^f - \overline{p}_n^f)
- \sum_i
T^{mf, UP}_{il} (\check{\overline{s}}^m_i, \check{\overline{s}}^f_l) (\overline{p}_i^m - \overline{p}_l^f )
 =
 \overline{q}_l^f  |\gamma_l|,
\end{split}
\end{equation}
where $l = 1,.., N^{f, H}$ and $i = 1,..., N^{m,H}$. For approximation by time similarly to the fine grid approximation, the IMPES scheme is used.

For the general multicontinuum model, we have
\begin{equation}
\label{tp-up2}
\begin{split}
\overline{\phi}^{\alpha}_i
\frac{\overline{s}^{\alpha}_i - \check{\overline{s}}^{\alpha}_i}{\tau} |K_i|
+ & \sum_{j}
T^{w,\alpha \alpha, UP}_{ij} (\check{\overline{s}}^{\alpha}_i, \check{\overline{s}}^{\alpha}_j)  (\overline{p}_i^{\alpha} - \overline{p}_j^{\alpha})
+  \sum_{\beta} \sum_l
T^{w, \alpha \beta, UP}_{il}(\check{\overline{s}}^{\alpha}_i, \check{\overline{s}}^{\beta}_l)
(\overline{p}_i^{\alpha} - \overline{p}_l^{\beta} )
=   (1 - \overline{f}^{w,\alpha}_i)
\overline{q}_i^{w,\alpha} |K_i|,  \\
& \sum_{j}
T^{\alpha \alpha, UP}_{ij} (\check{\overline{s}}^{\alpha}_i, \check{\overline{s}}^{\alpha}_j)  (\overline{p}_i^{\alpha} - \overline{p}_j^{\alpha})
+  \sum_{\beta} \sum_l
T^{\alpha \beta, UP}_{il}(\check{\overline{s}}^{\alpha}_i, \check{\overline{s}}^{\beta}_l) (\overline{p}_i^{\alpha} - \overline{p}_l^{\beta} )
=
\overline{q}_i^{\alpha} |K_i|, \\
\end{split}
\end{equation}
where
\begin{equation}
\label{tp-up3}
T_{ij}^{\alpha \beta, UP}(\overline{s}^{\alpha \beta}_{ij}) =
\lambda(\overline{s}^{\alpha}_i, \overline{s}^{\beta}_j)
W_{ij}^{\alpha \beta, UP}, \quad
T_{ij}^{w, \alpha \beta, UP}(\overline{s}^{\alpha}_i, \overline{s}^{\alpha}_j) =
\lambda^w (\overline{s}^{\alpha \beta}_{ij})
W_{ij}^{\alpha \beta, UP},  \quad \alpha, \beta = m,f
\end{equation}
with upwind scheme approximation of $\lambda^w$ and $W^{\alpha \beta, UP}_{ij}$ is the precalculated effective transmissibilities that is similar to the previous problem and based on the single phase upscaling.


The choice of boundary conditions have a strong impact on the accuracy of  results.
In the presented standard upscaling method, coarse grid parameters are obtained independently to  global problem solution information.
More accurate approaches can be based on the information about the fine scale flow in the local domains up to fine grid resolution and without variable separation of nonlinear coefficients.
For example, an interpolated global coarse grid solution is used for performing accurate construction of the upscaled transmissibilities  in \cite{chen2003coupled}, which involve iterations between global coarse grid model and local fine grid calculations with updating of the upscaled transmissibilities. The local-global upscaling method requires extra computations than existing classic upscaling procedures.

In this work, the construction of the accurate upscaled transmissibilities for the coarse grid approximation  is also based on the information about global solution (nonlinear transmissibilities).
Moreover, the presented method is based on the machine learning procedure for fast prediction of the nonlinear transmissibilities, where we construct neural network that learn dependencies between the coarse grid quantities on the oversampled local domains and upscaled transmissibilities. We use a convolutional neural network and GPU training process to construct a machine learning algorithm.

\section{Machine learning for nonlinear nonlocal upscaled transmissibilities}

We consider a machine learning approach for prediction of the upscaled nonlinear nonlocal transmissibilities for accurate and fast coarse grid approximation.
We have following main steps:
\begin{enumerate}
\item Generate dataset to train, validate and test of the neural network.
\item Neural networks training, validation and testing.
\item Calculation of the nonlinear upscaled transmissibilities on the fly using constructed neural networks during coarse system construction, fast and accurate solution of the upscaled system.
\end{enumerate}

For construction of the datasets, we perform local or global calculations of the coarse grid quantities \cite{chen2003coupled}. In local approach, the upscaled transmissibilities are calculated on the local domain corresponding to the target face, where the fine-scale solution information is used to set boundary conditions. Global approach uses a global fine-scale solution for the determination of coarse scale parameters.
For training of the neural networks, we use a family of problem solutions for different input conditions (snapshots). Note that, we should have many snapshots to capture all input condition variations because the accuracy of the machine learning method depends on snapshot space that is as train dataset.
Next, we consider dataset generation and
network construction  in detail.

\subsection{Dataset}

The most accurate case can be based on the fine grid solution, at the same time for upscaled model, we would like to use only coarse-grid information. For possible applicability of this, we construct a novel coarse grid model, using machine learning algorithms and construct neural network that learn dependency between coarse grid functions $\overline{p}^{\alpha}$ in local domains and upscaled nonlinear transmissibilities.

For constructing accurate neural network for prediction of the transmissibilities, we should train network on the highly accurate dataset.
One of the most accurate approach for calculating upscaled transmissibilities based
on the direct calculation from the fine scale solution.
We use following coarse grid approximation (similar to previous section)
\begin{itemize}
\item Unsaturated flow problem (nonlinear flow):
\begin{equation}
\label{ml-uns}
\begin{split}
\overline{c}^{\alpha}_i
\frac{ \overline{p}_i^{\alpha} - \check{\overline{p}}_i^{\alpha} }{\tau} |K_i|
+ \sum_{j}
{T}_{ij}^{\alpha \alpha, NL} (\overline{p}_i^{\alpha} - \overline{p}_j^{\alpha})
+ \sum_{\beta} \sum_l
{T}^{\alpha \beta, NL}_{il} (\overline{p}_i^{\alpha} - \overline{p}_l^{\beta} )
= \overline{q}^{\alpha}_i |K_i|,
\end{split}
\end{equation}
with nonlinear upscaled transmissibilities
\begin{equation}
\label{ml-uns2}
{T}_{ij}^{\alpha \beta, NL}(x, p^{\alpha}, p^{\beta}, s^{\alpha}, s^{\beta}) =
\frac{ \sum_{r,n}
T^{\alpha \beta}_{rn}(s^{\alpha}, s^{\beta}) (p^{\alpha}_r - p^{\beta}_n) }{\overline{p}^{\alpha}_i - \overline{p}^{\beta}_j}.
\end{equation}
\item Two-phase flow problem (nonlinear transport and flow):
\begin{equation}
\label{ml-tp}
\begin{split}
\overline{\phi}^{\alpha}_i
\frac{\overline{s}^{\alpha}_i - \check{\overline{s}}^{\alpha}_i}{\tau} |K_i|
+ & \sum_{j}
T^{w,\alpha \alpha, NL}_{ij} (\overline{p}_i^{\alpha} - \overline{p}_j^{\alpha})
+  \sum_{\beta} \sum_l
T^{w, \alpha \beta, NL}_{il} (\overline{p}_i^{\alpha} - \overline{p}_l^{\beta} )
=   (1 - \overline{f}^{w,\alpha}_i)
\overline{q}_i^{w,\alpha} |K_i|,  \\
& \sum_{j}
T^{\alpha \alpha, NL}_{ij}  (\overline{p}_i^{\alpha} - \overline{p}_j^{\alpha})
+  \sum_{\beta}  \sum_l
T^{\alpha \beta, NL}_{il}  (\overline{p}_i^{\alpha} - \overline{p}_l^{\beta} )
=
\overline{q}_i^{\alpha} |K_i|, \\
\end{split}
\end{equation}
with nonlinear upscaled transmissibilities
\begin{equation}
\label{ml-tp2}
\begin{split}
{T}_{ij}^{w, \alpha \beta, NL}(x, p^{\alpha}, p^{\beta}, s^{\alpha}, s^{\beta})  =
\frac{ \sum_{r,n}
T^{w, \alpha \beta}_{rn}(s^{\alpha}, s^{\beta}) (p^{\alpha}_r - p^{\beta}_n) }{\overline{p}^{\alpha}_i - \overline{p}^{\beta}_j},
\\
{T}_{ij}^{\alpha \beta, NL}(x, p^{\alpha}, p^{\beta}, s^{\alpha}, s^{\beta}) =
\frac{ \sum_{r,n}
T^{\alpha \beta}_{rn}(s^{\alpha}, s^{\beta}) (p^{\alpha}_r - p^{\beta}_n) }{\overline{p}^{\alpha}_i - \overline{p}^{\beta}_j}.
\end{split}
\end{equation}
\end{itemize}

Because ${T}^{w, \alpha \beta, NL}$ and ${T}^{\alpha \beta, NL}$ ($\alpha, \beta = m,f$) are nonlinear and depend on the fine grid solutions $p^{\alpha}, p^{\beta}, s^{\alpha}, s^{\beta}$,
we cannot directly use such transmissibilities on the coarse grid model. For possible applicability of this, we will use a machine learning algorithms and construct a neural network that learn dependence between coarse grid functions $\overline{p}^{\alpha}, \overline{p}^{\beta}, \overline{s}^{\alpha}, \overline{s}^{\beta}$ in local domains (oversampled) and upscaled nonlinear transmissibilities.

\begin{figure}[h!]
\centering
\includegraphics[width=0.5\linewidth]{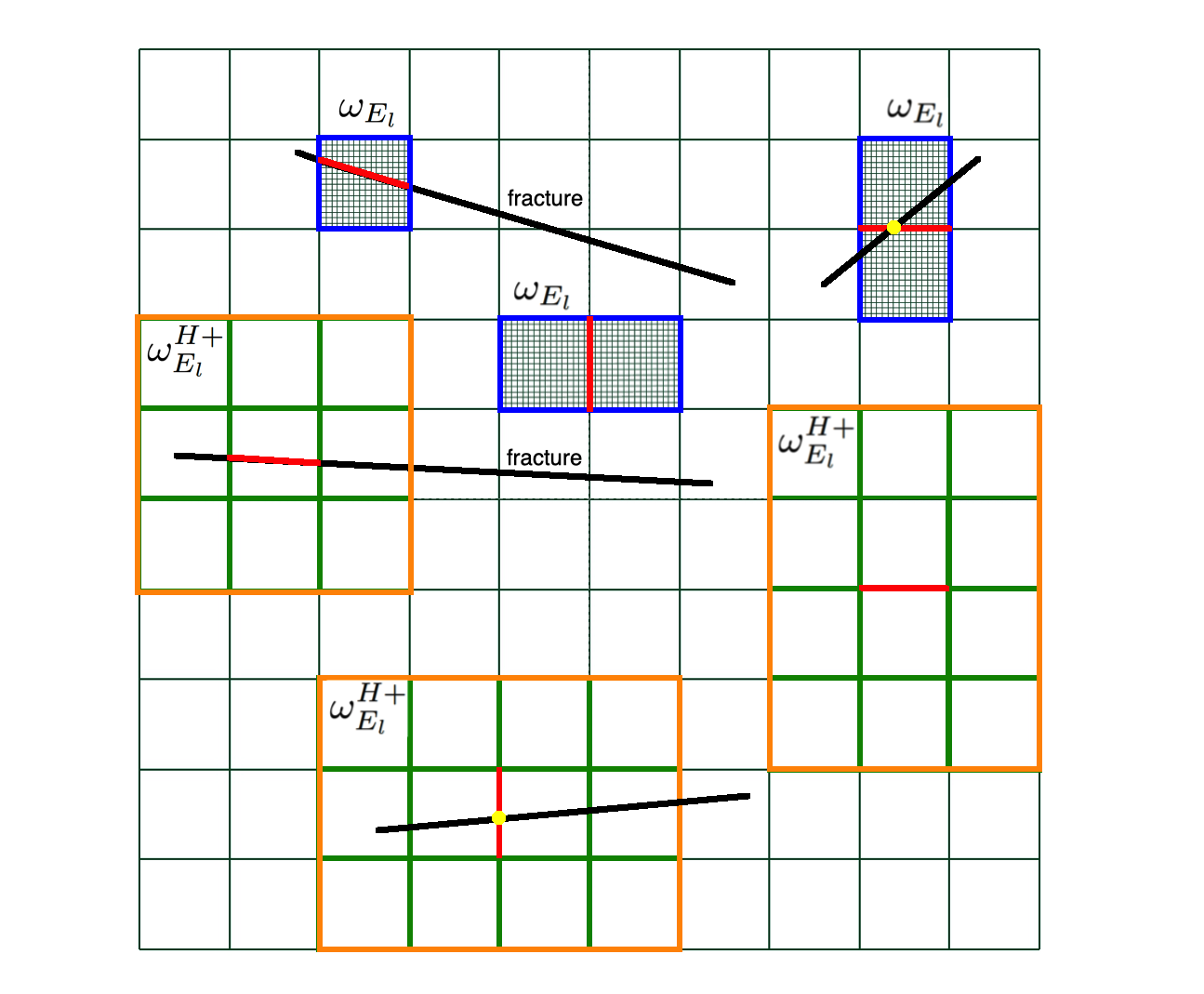}
\caption{Coarse grid and local domains illustration.
Four types (horizontal and vertical matrix-matrix flow, fracture-matrix flow and fracture-fracture flow) of local domains $\omega_{E_l}$ ($\partial \omega_{E_l}$ -- blue color, fine grid resolution -- black color, coarse edge or matrix-fracture interface -- red color and fracture-fracture connection -- yellow point).
Four types (horizontal and vertical matrix-matrix flow, fracture-matrix flow and fracture-fracture flow) of oversampled local domains $\omega^{H+}_{E_l}$ ($\partial \omega^{H+}_{E_l}$ -- orange color, coarse grid resolution -- green color).  }
\label{fig:omega}
\end{figure}

Let $E_l$ is the interface, where we define upscaled transmissibility, and  $X_l$ and $Y_l$ are the input data and output data for machine learning algorithm and
\[
\text{ Dataset: }
\{ (X_l, Y_l), \, l = 1,...,L \}.
\]
For constructing neural network for upscaled transmissibilities, based on the \eqref{ml-uns2} and \eqref{ml-tp2}, we use a nonlocal upscaled transmissibilities ${T}_l^{w, \alpha \beta, NL}$ and ${T}_l^{\alpha \beta, NL}$ ($\alpha, \beta = m,f$) as output data $Y_l$.
Input data $X_l$ is contains information about fine scale permeabilities, fracture position in local domain $\omega_{E_l}$, coarse grid functions $\overline{p}$ and $\overline{s}$ in the oversampled local domains.
For this purpose, we use a local multi-input data for training neural network
\begin{equation}
\label{ml-xy1}
X_l = (X_l^k, X_l^f, X_{l+}^{\overline{p}^{\alpha}}, X_{l+}^{\overline{p}^{\beta}}) \text{ for flow} \quad \text{ and } \quad
X_l = (X_l^k, X_l^f, X_{l+}^{\overline{p}^{\alpha}}, X_{l+}^{\overline{p}^{\beta}}, X_{l+}^{\overline{s}^{\alpha}}, X_{l+}^{\overline{s}^{\beta}}) \text{ for transport and flow},
\end{equation}
where
$X_l^k$ and $X_l^f$ are the local heterogeneous permeabilities and local fracture position markers in local domain $\omega_{E_l}$;
$X_{l+}^{\overline{p}^{\alpha}}$ and $X_{l+}^{\overline{s}^{\alpha}}$ are the coarse grid nonlocal mean values for pressure and saturation for continuum $\alpha$ in oversampled local domain $\omega^{H+}_{E_l}$. Each of the input fields is represented as two-dimensional array for two-dimensional problem. The scale of each array in dataset is re-scaled to fall within the range $0$ to $1$.

In Figure \ref{fig:omega}, we present an illustration of the local domains $\omega_{E_l}$ and $\omega^{H+}_{E_l}$, $l = 1,...,N_E$ ($N_E$ -- number local domains).
Local domain $\omega_{E_l}$ is the domain for edge $E_l$ up to fine grid resolution that is similar to the classic (single phase) upscaling presented in previous section. In local domain $\omega_{E_l}$, we define $X_l^k \in \omega_{E_l}$ and $X_l^f \in \omega_{E_l}$.
Oversampled local domain $\omega^{H+}_{E_l}$ is the domain around $E_l$ up to coarse grid resolution, where we define $X_{l+}^{\overline{p}^{\alpha}} \in \omega^{H+}_{E_l}$ and $X_{l+}^{\overline{s}^{\alpha}} \in \omega^{H+}_{E_l}$.
To ensure same size and structure of the input data, we divide all local data into four types: matrix-matrix flow through horizontal edge ($T_l^{mm, NL}$), matrix-matrix flow through vertical edge ($T_l^{mm, NL}$, fracture -matrix flow  for $T_l^{mf, NL}$ and fracture-fracture flow for $T_l^{ff, NL}$.

The output  is the normalized array of the upscaled transmissibilities
\[
Y_l = (T_l^{\alpha \beta, NL})
\text{ for flow} \quad \text{ and } \quad
Y_l = (T_l^{\alpha \beta, NL}, T_l^{w, \alpha \beta, NL}) \text{ for transport and flow},
\]
for edge $E_l$.
Dataset is divided into train, validation and test sets with sizes $N_{train}$, $N_{val}$ and $N_{test}$ ($N = N_{train} + N_{val} + N_{test}$). For each type of local domain as a test set,
we take 50 \% of data, another 50 \%  divided between train and validation
set in 80/20 proportion.

We use dataset $(X_l, Y_l)$ for training of the neural network. To ensure a good learning rate and for obtaining a wide coverage of data, we generate several solution snapshots by varying of the source term in  global fine grid model. Another approach is related to the localization of the dataset generation, where we can use local domains calculations for fine grid solutions and calculations of the $X_l$ and $Y_l$. Note that, this machine learning approach for the learning of the nonlocal nonlinear upscaled transmissibilities can be also applied for the linear problems and has a recap with nonlinear finite volume methods, where transmissibilities are also depends on solution.

\subsection{Network}

\begin{figure}[h!]
\centering
\includegraphics[width=0.5\linewidth]{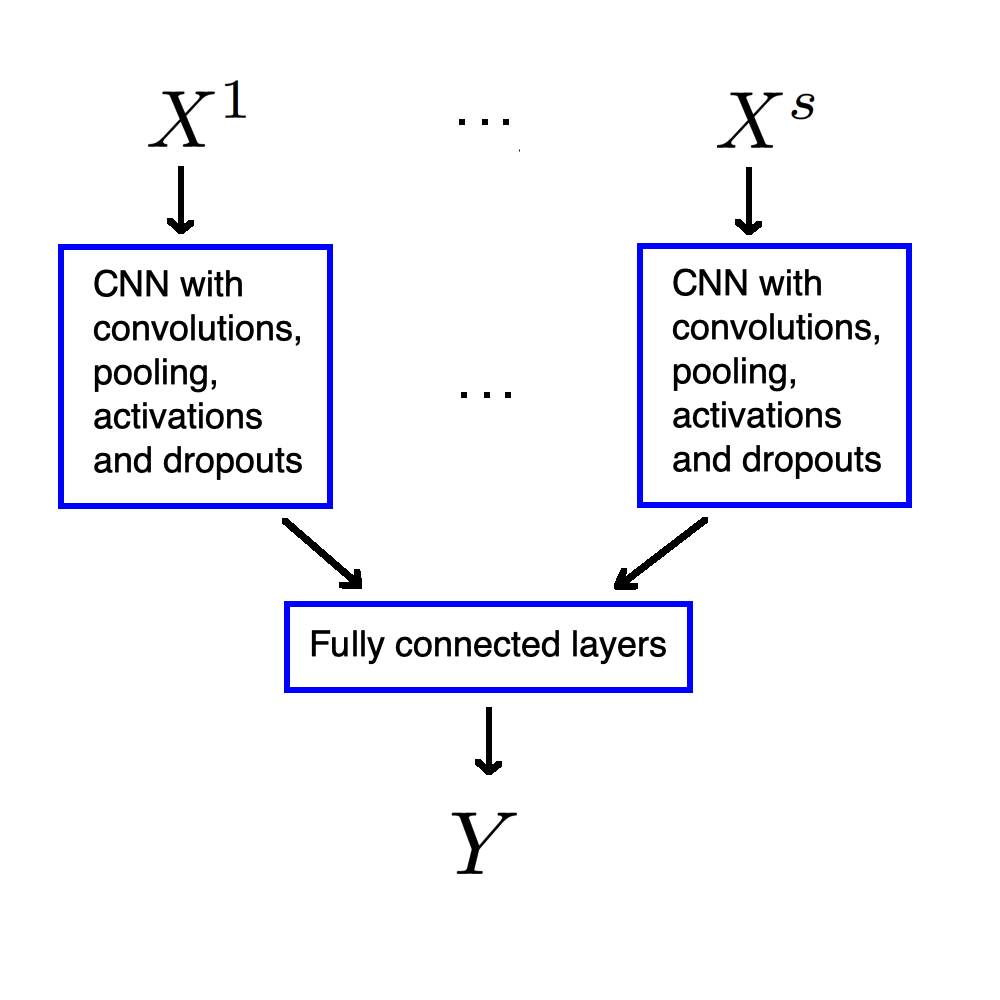}
\caption{Illustration of the multi-input deep neural network for prediction of the nonlinear nonlocal upscaled transmissibilities }
\label{fig:ml}
\end{figure}

In machine learning algorithm, we use a multi-input deep neural network (convolutional neural network).
Let
\[
\text{ Dataset: }
\{ (X_l, Y_l), \, l = 1,...,L \}
\]
where $X_l = (X_l^1, ..., X_l^s)$ ($s$ is the number of the input data for $E_l$, see \eqref{ml-xy1}). Each input data $X_l^i$ is defined in $\omega^i_{E_l}$ and represented as two-dimensional array for two-dimensional problems.
The architecture of the multi-input deep neural network for prediction of the nonlinear nonlocal upscaled transmissibilities is presented in Figure \ref{fig:ml}.
For each input data $X^i$, we use a convolutional neural network \cite{lecun2015deep, krizhevsky2012imagenet}. Several convolutional and pooling layers with rectified linear units activation layer are stacked with a several fully-connected layers with dropout. Several layers of convolutions and pooling are alternated in order to detect higher order features for better accuracy of the method.
After performing convolutions, pooling, activation and dropout layers for each $X^i$ ($i = 1,...,s$), we add a fully connected layers, where we compose all outputs on CNN together.
By a training process, a machine learning algorithm solve the optimization problem to find model weights that best describe the train set by minimization of the loss function.

We train a convolutional neural network by a dataset of local multi-input data ($X_l$) and upscaled transmissibilities ($Y_l$). As a loss function, we use the mean square error (MSE)
\[
Loss_{train} =
\frac{1}{N_{train}} \sum_{l=1}^{N_{train}} |Y_l - F(X_l) |^2.
\]
For solution of the minimization problem, we use gradient-based optimizer Adam \cite{kingma2014adam}. Implementation of the machine learning method is based on the open source library Keras \cite{keras} with TensorFlow backend \cite{tensorflow} and performed on the GPU.
Constructed machine learning algorithm will efficiently determine dependence between coarse grid functions in local domains and upscaled transmissibilities.

\section{Numerical result}

In this section, we present numerical results for the proposed method. We consider following model problems in fractured and heterogeneous porous media:
\begin{itemize}
\item[] \textit{Test 1}: Nonlinear flow problem (unsaturated flow problem)
\item[] \textit{Test 2}: Nonlinear transport and flow problem (two-phase flow problem)
\end{itemize}

\begin{figure}[h!]
\centering
\includegraphics[width=0.32\linewidth]{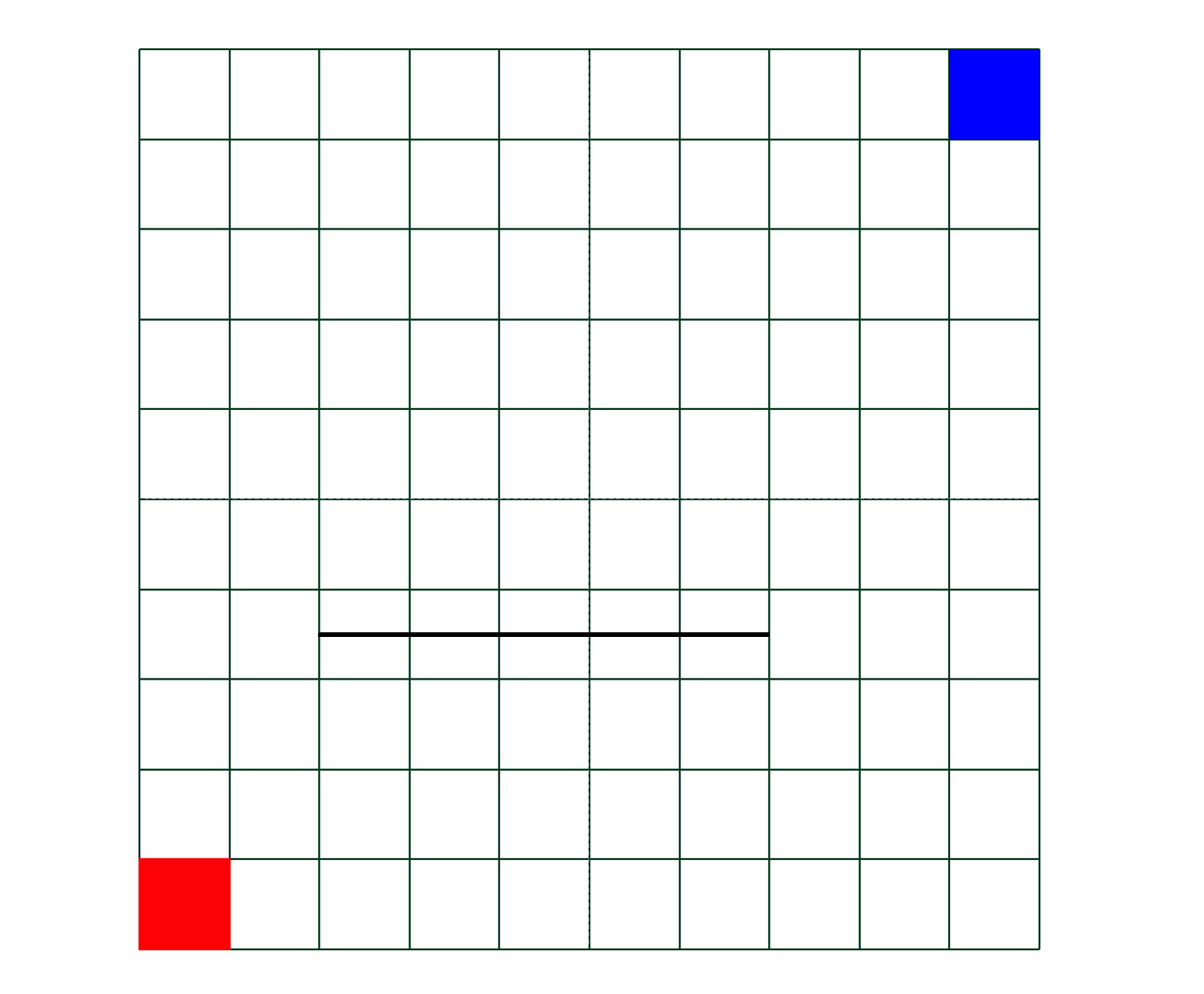}
\includegraphics[width=0.3\linewidth]{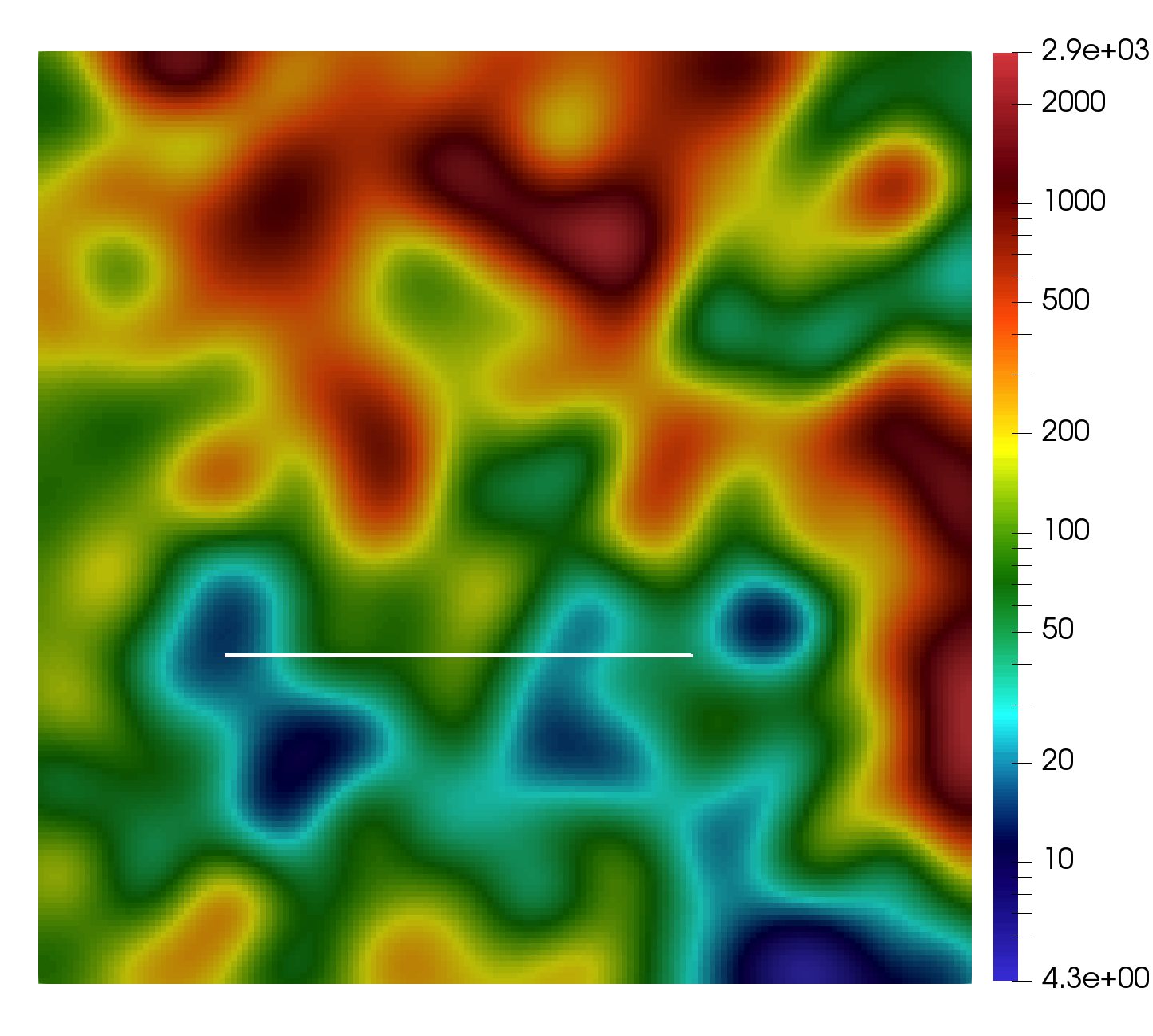}
\includegraphics[width=0.3\linewidth]{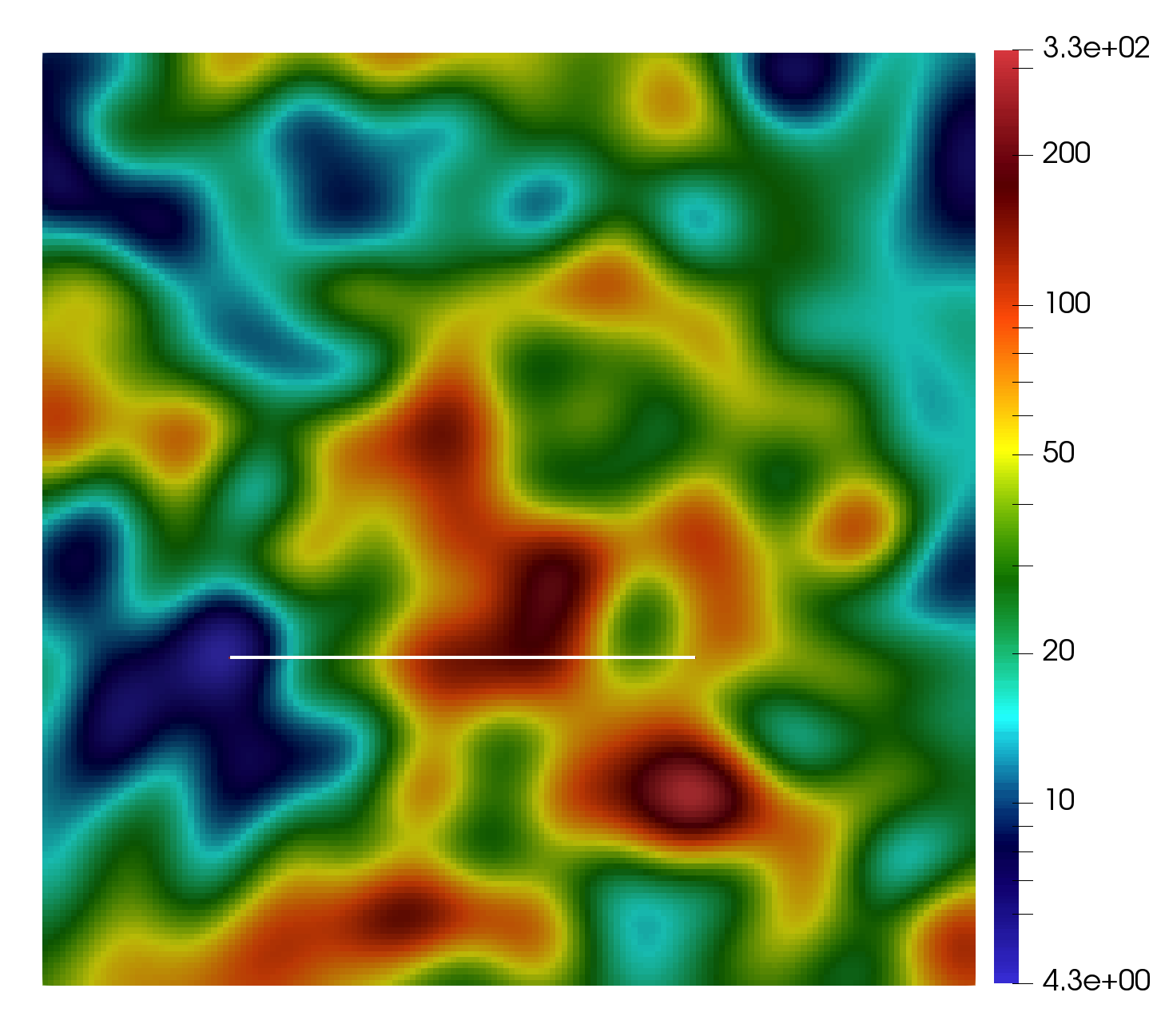}
\caption{Coarse mesh with source term and fracture positions (left).  Heterogeneous porous matrix permeability in $\Omega$ for \textit{Test 1} and \textit{Test 2}}
\label{fig:kx}
\end{figure}

We solve model problem in $\Omega = [0, 1] \times  [0, 1]$ with no flux boundary conditions.  We use $10 \times 10$ coarse grid.
Location of the source terms and fracture position are depicted in  Figure \ref{fig:kx}. In Figure \ref{fig:kx}, we show a heterogeneous porous matrix permeability for both test problems.
The numerical calculations of the effective properties has been implemented with the open-source finite element software PETSc and FEniCS \cite{fenics, logg2012automated, balay2004petsc}.

To measure difference between reference solution and coarse grid solution, we compute relative $L_2$ error
\[
e(\overline{u}) = \sqrt{
\frac{\sum_{i = 1}^{N^H} (\overline{u}^{fine}_i - \overline{u}_i)^2 }{\sum_{i = 1}^{N^H}  (\overline{u}^{fine}_i)^2 }
},
\]
where $u = p, s$,  $\overline{u}^{fine}$ is the reference solution (mean value on coarse grid of the fine grid solution) and $\overline{u}$ is the solution on the coarse grid.

For each test problem, we present results of the fine scale solution, for  upscaling technique presented in Section 3 and new method from Section 4.
Computational algorithm for single-phase  upscaling method with $T^{\alpha \beta, UP}$ (Section 3):
\begin{enumerate}
\item Loading of the precalculated effective transmissibilities $W^{\alpha \beta, UP}$.
\item Solution of the multicontinuum model:
\begin{itemize}
\item[] \textit{Test 1}: Nonlinear flow problem with
\[
{T}^{\alpha \beta, UP}_{ij}(\overline{p}_i^{\alpha}, \overline{p}_j^{\beta})
=
k_r(\overline{p}_{ij}^{\alpha \beta})
W^{\alpha \beta, UP}_{ij}, \quad \alpha = m, f
\]
\item[] \textit{Test 2}: Nonlinear transport and flow problem with
\[
T_{ij}^{\alpha \beta, UP}(\overline{s}^{\alpha}_i, \overline{s}^{\beta}_j) =
\lambda(\overline{s}_{ij}^{\alpha \beta})
W_{ij}^{\alpha \beta, UP}, \]
\[
T_{ij}^{w, \alpha \beta, UP}(\overline{s}^{\alpha}_i, \overline{s}^{\alpha}_j) =
\lambda^w (\overline{s}_{ij}^{\alpha \beta})
W_{ij}^{\alpha \beta, UP},  \quad \alpha, \beta = m,f
\]
with upwind  approximation of $\lambda^w$ on the coarse grid.
\end{itemize}
\end{enumerate}
For the new nonlocal nonlinear machine learning technique with $T^{\alpha \beta, NL}$ (Section 4), we have:
\begin{enumerate}
\item Loading of the machine learning models, $NN_i$ ($i=1,2,...$).
\item Solution of the multicontinuum model:
\begin{itemize}
\item[] \textit{Test 1}: Nonlinear flow problem with
\[
{T}_{ij}^{\alpha \beta, NL}  =
\begin{cases}
{T}_{ij}^{\alpha \beta, ML}(x, \overline{p}^{\alpha}, \overline{p}^{\beta}, \overline{p}^{\alpha}, \overline{p}^{\beta}), & \text{ if } |\overline{p}^{\alpha} - \overline{p}^{\beta}| > \varepsilon, \\
{T}_{ij}^{\alpha \beta, UP}, & else
\end{cases}.
\]
where ${T}_{ij}^{\alpha \beta, ML}$ is the value predicted using machine learning algorithm.
\item[] \textit{Test 2}: Nonlinear transport and flow problem with
\[
{T}_{ij}^{\alpha \beta, NL}  =
\begin{cases}
{T}_{ij}^{\alpha \beta, ML}(x, \overline{p}^{\alpha}, \overline{p}^{\beta}, \overline{p}^{\alpha}, \overline{p}^{\beta}), & \text{ if } |\overline{p}^{\alpha} - \overline{p}^{\beta}| > \varepsilon, \\
{T}_{ij}^{\alpha \beta, UP}, & else
\end{cases},
\]\[
{T}_{ij}^{w, \alpha \beta, NL}  =
\begin{cases}
{T}_{ij}^{w, \alpha \beta, ML}(x, \overline{p}^{\alpha}, \overline{p}^{\beta}, \overline{p}^{\alpha}, \overline{p}^{\beta}), & \text{ if } |\overline{p}^{\alpha} - \overline{p}^{\beta}| > \varepsilon, \,   |\overline{s}^{\alpha} - \overline{s}^{\beta}|  > \varepsilon_s\\
{T}_{ij}^{w, \alpha \beta, UP}, & else
\end{cases}.
\]
where ${T}_{ij}^{\alpha \beta, ML}$ and ${T}_{ij}^{w, \alpha \beta, ML}$ is the value predicted using machine learning algorithm.
\end{itemize}
\end{enumerate}
Note that, the loss of positivity of upscaled transmissibilities can happen, and we use a threshold value $\varepsilon$ for the pressure difference to guarantee a good values of the coarse grid parameters, where linear upscaling is used for the faces with small pressure difference.
Moreover, we used predicted transmissibilities adaptively with parameter $\varepsilon_s$ in the coarse grid model with machine learning approach.

We will show results of the learning process of deep neural network for nonlocal nonlinear upscaled transmissibilities  and calculate errors for a given datasets.
Finally, we consider a coarse grid solution of the problem, where nonlocal nonlinear upscaled transmissibilities are calculated using constructed machine learning method. 
Finally, we discuss the computational time of the neural networks construction and solution of the coarse grid system using classic upscaling and machine learning approaches.
We divide calculation on the offline and online stages. On the online stage, we train neural network on the GPU by a given train and validation datasets.
On the offline stage, we have two steps: loading of the preconstructed neural network and prediction of the upscaled coarse grid transmissibilities on each time iteration or/and nonlinear iteration.

\subsection{Nonlinear flow problem}

We consider the solution of the nonlinear equation in fractured and heterogeneous porous media. We set source terms $f^{\pm} = \pm q$, $q = 10^5$. For the nonlinear coefficient, we use $k^{\alpha \beta}(x, u) = k_s(x) k_r(u)$ with $k_r(u) = \exp(-a |u|)$, $a = 0.1$ ($\alpha, \beta = m,f$).
 In Figure \ref{fig:kx} (second column),  we show a heterogeneous porous matrix permeability $k^m_s(x)$ and fracture position. We set $c^m = 1$, $c^f = 0$, $k_s^f = 10^6$ and $T_{max} = 10^{-3}$ with 20 time steps. Coarse grid is $10 \times 10$ and fine grid is $640 \times 640$ for domain $\Omega$.

\begin{figure}[h!]
\centering
\includegraphics[width=0.32\linewidth]{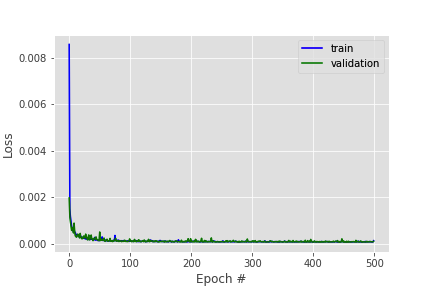}
\includegraphics[width=0.32\linewidth]{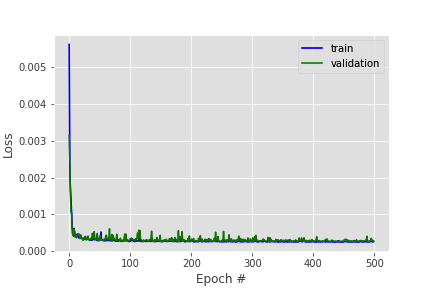}
\includegraphics[width=0.32\linewidth]{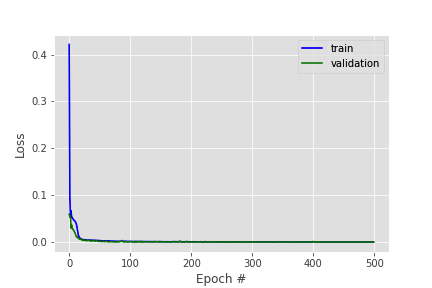}
\caption{Learning process  for \textit{Test 1} (nonlinear flow). Loss function vs epoch.
Left: $NN_1$ for vertical ${T}^{mm, NL}$.
Middle: $NN_2$ for horizontal ${T}^{mm, NL}$.
Right: $NN_3$ for  ${T}^{mf, NL}$ }
\label{fig:ml1}
\end{figure}

\begin{figure}[h!]
\centering
\includegraphics[width=0.32\linewidth]{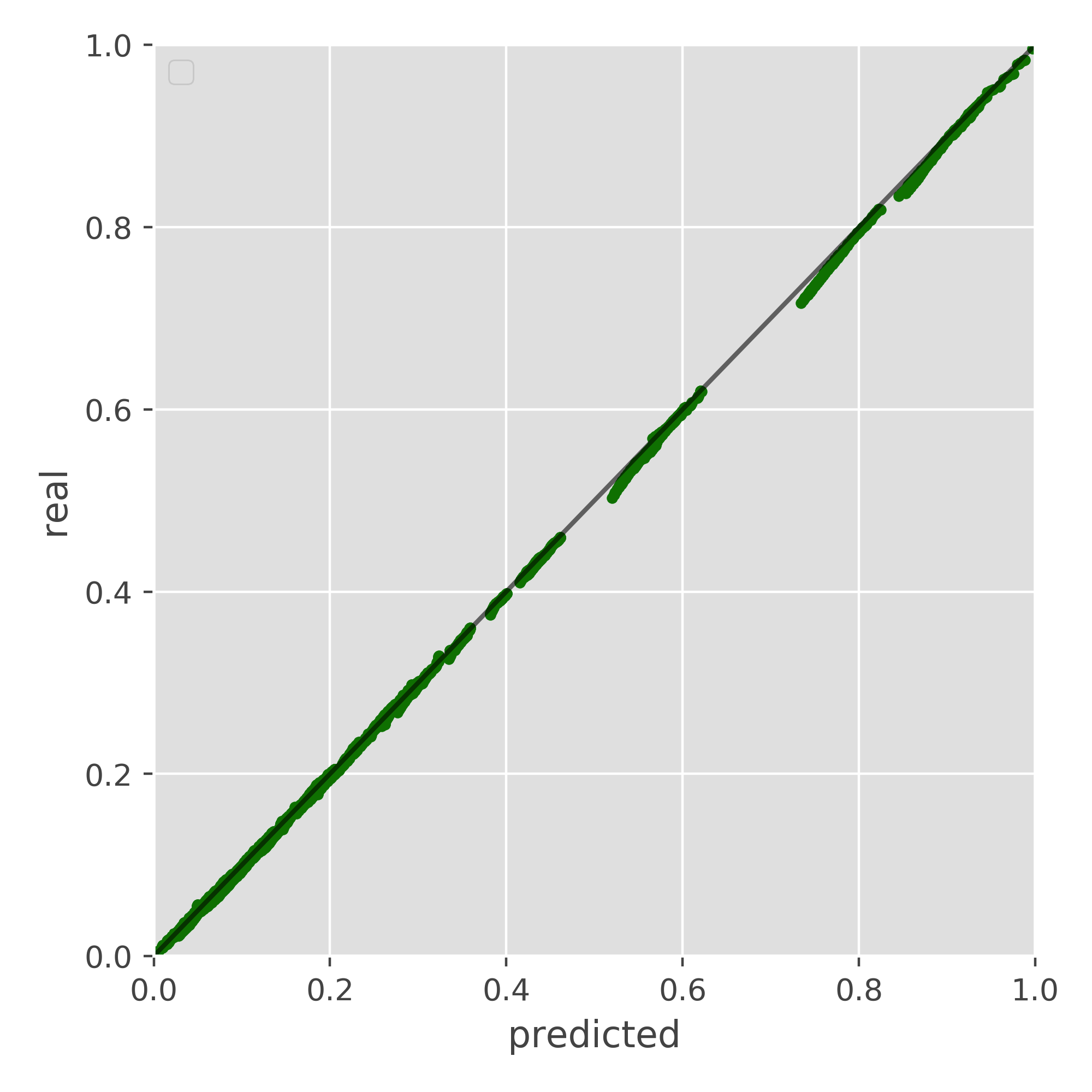}
\includegraphics[width=0.32\linewidth]{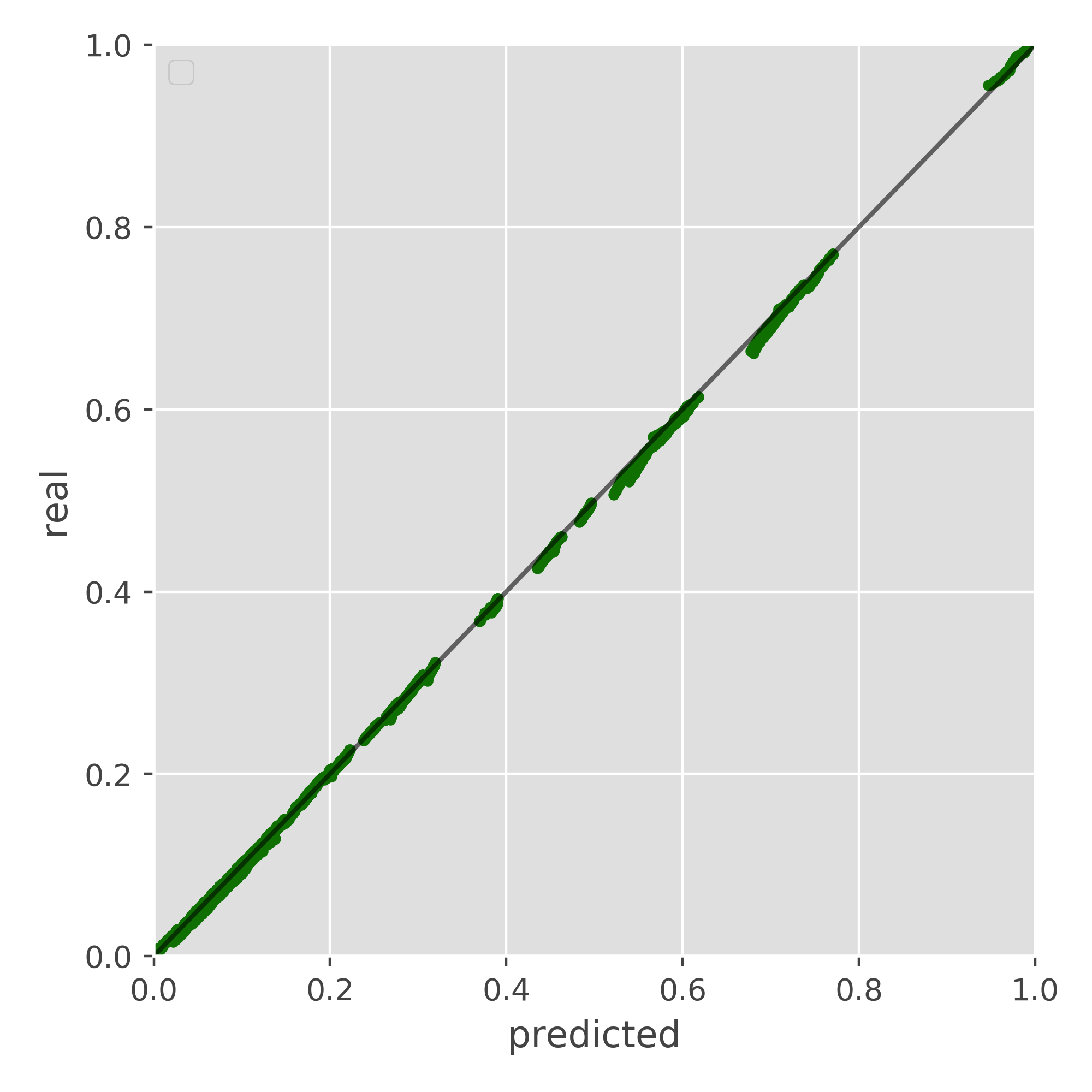}
\includegraphics[width=0.32\linewidth]{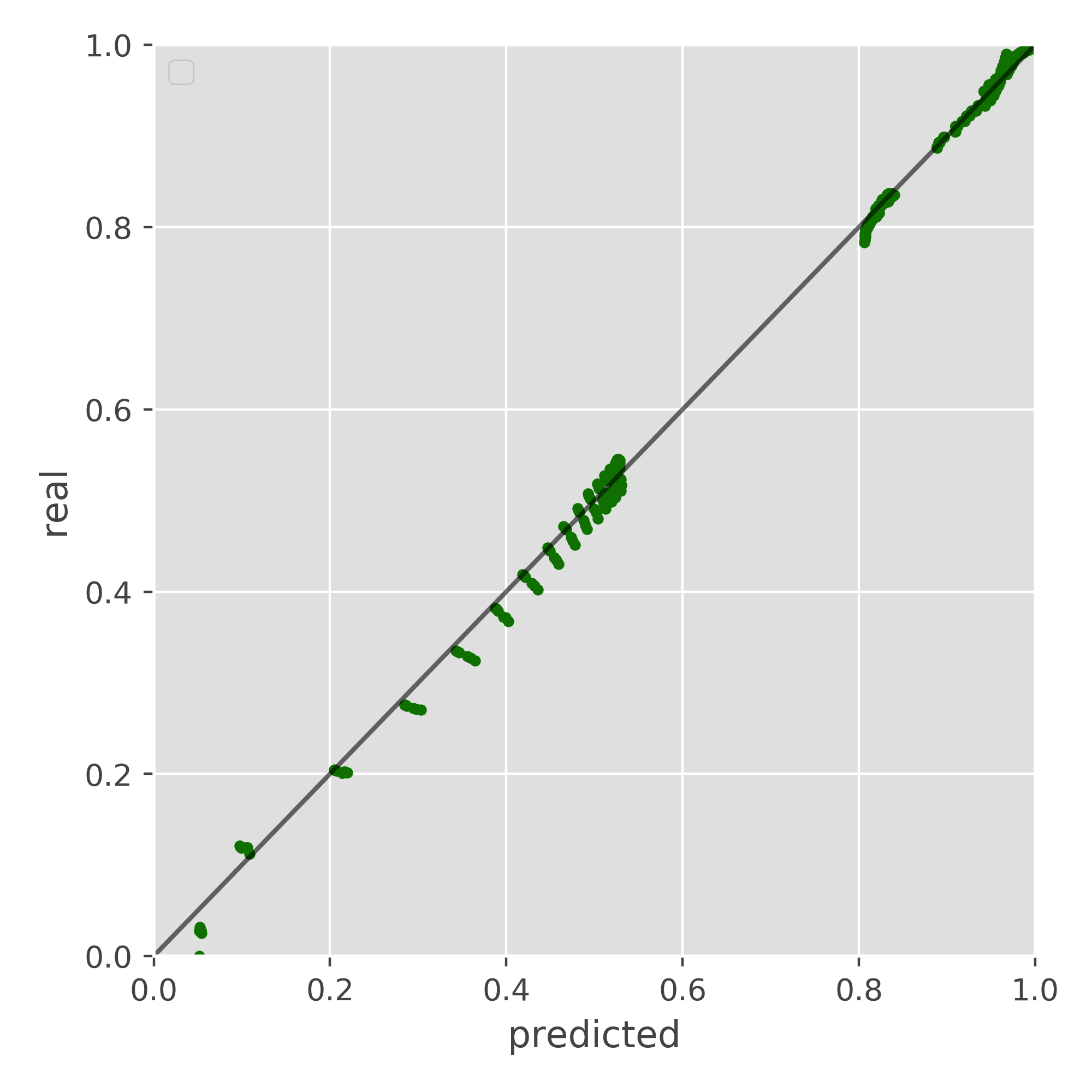}\\
\includegraphics[width=0.32\linewidth]{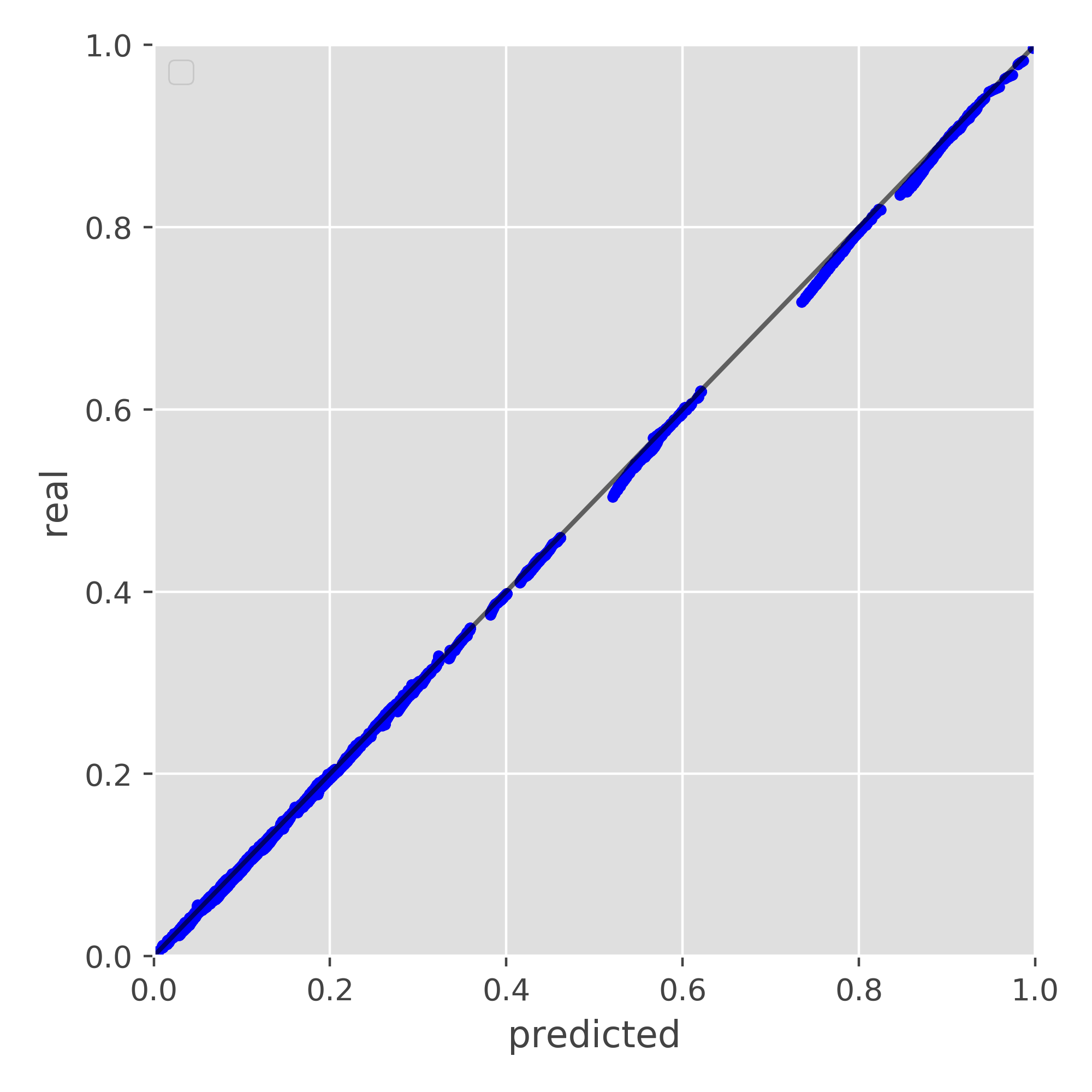}
\includegraphics[width=0.32\linewidth]{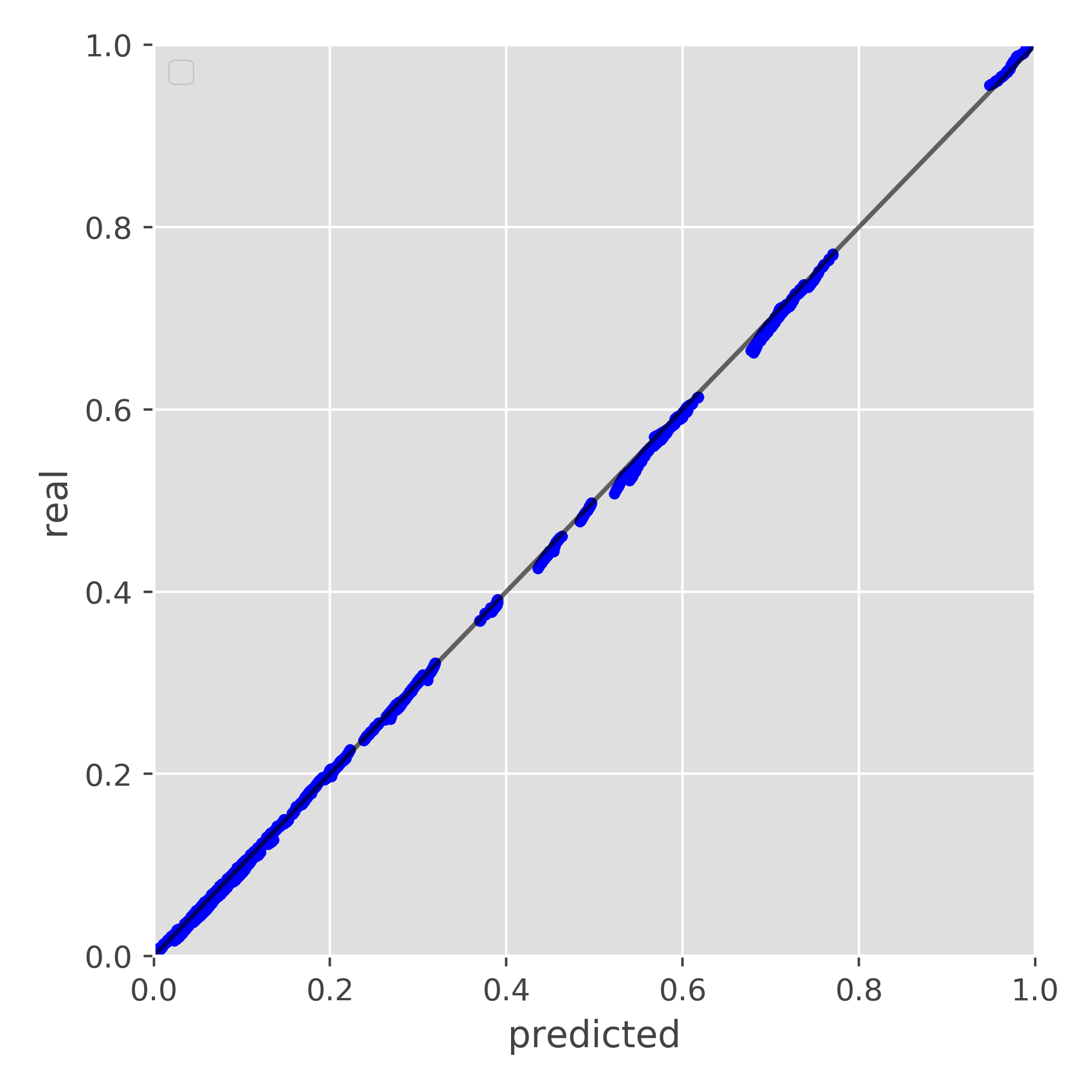}
\includegraphics[width=0.32\linewidth]{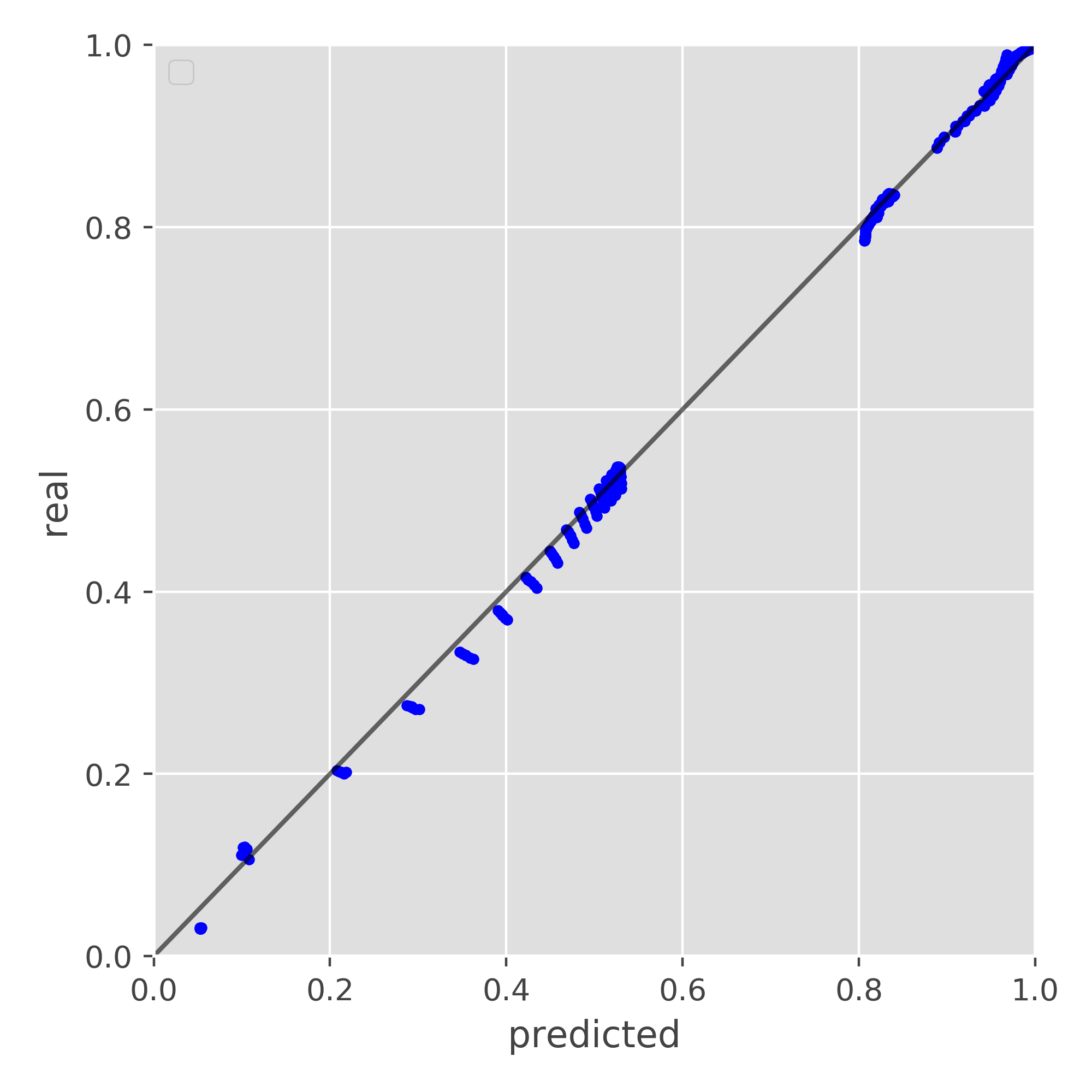}
\caption{Learning performance  for \textit{Test 1} (nonlinear flow). Parity plots comparing preference property values against predictions made using machine learning.
First row: train and validation dataset (green color). Second row: test dataset (blue color).
Left: $NN_1$ for vertical ${T}^{mm, NL}$.
Middle: $NN_2$ for horizontal ${T}^{mm, NL}$.
Right: $NN_3$ for  ${T}^{mf, NL}$ }
\label{fig:ml2}
\end{figure}

\begin{table}[h!]
\begin{center}
\begin{tabular}{ | c | c c c | }
\hline
& MSE & RMSE  (\%)  & MAE (\%) \\
\hline
\multicolumn{4}{|c|}{Train set (global data)} \\
\hline
$NN_1$ 	& 0.012 & 1.101 & 1.072 \\
$NN_2$ 	& 0.016 & 1.273 & 1.272 \\
$NN_3$ 	& 0.013 & 1.156 & 0.838 \\
\hline
\multicolumn{4}{|c|}{Test set (global data)} \\
\hline
$NN_1$ & 0.012 & 1.104 & 1.085 \\
$NN_2$ & 0.016 & 1.286 & 1.280 \\
$NN_3$ & 0.011 & 1.050 & 0.774 \\
\hline
\multicolumn{4}{|c|}{Train set (local data)} \\
\hline
$NN^l_1$ 	& 0.081 & 2.861 & 2.060 \\
$NN^l_2$ 	& 0.014 & 1.223 & 1.229 \\
$NN^l_3$ 	& 0.042 & 2.058 & 1.791 \\
\hline
\end{tabular}
\end{center}
\caption{Learning performance of machine learning algorithm for \textit{Test 1} (nonlinear flow). Errors for train and test sets }
\label{tab:ml}
\end{table}

We present results for the machine learning algorithm and calculate errors for train and test datasets. For the training of the neural networks, we investigate two datasets: local and global. For the global dataset, we extract local information from the fine grid calculations on the global domain $\Omega$. For the local dataset, we calculate each data by solution of the local problem up to fine grid resolution with different boundary conditions for generation of the possible set of solutions (snapshots).
We use six random values of the  source term  to generate datasets ($N_r = 6$).
We train three neural networks for each type of transmissibility: $NN_1$ for horizontal coarse edges for matrix-matrix flow, $NN_2$ for vertical coarse edges s for matrix-matrix flow and $NN_3$ for matrix - fracture flow. For $10 \times 10$ coarse mesh, we have $N_E = 90$ horizontal and $N_E = 90$ vertical coarse edges (without boundary edges due to no flux boundary conditions), furthermore, we have $N_E = 5$ coarse cells with fracture.
Therefore, the train dataset for  neural network contains $N = N_r \cdot N_E \cdot N_t$ samples for learning process, where $N_t$ is the number of time steps. We have $N = 10800$ for $NN_1$ and $NN_2$; and $N = 600$ for $NN_3$. Each sample $X_l$ contains information about heterogeneous permeability and fracture position up to fine grid resolution in local domain, coarse grid mean value of the solution in oversampled local domain
\[
X_l = (X_l^k, X_l^f, X_{l+}^{\overline{p}^{m}}).
\]
Each dataset is divided into training and validation sets with $80:20$ ratio. For testing, we calculate another six solution snapshots.

For calculations, we use 500 epochs with a batch size $N_b = 90$ and Adam optimizer with learning rate $\epsilon = 0.001$.
For accelerating of the training process of the multi-input CNN, we use GPU.
We use $3 \times 3$ convolutions and $2 \times 2$ maxpooling layers with RELU activation for $X^k$ and $X^f$, and $3 \times 3$ convolutions with RELU activation for $X^{\overline{p}^{m}}$. For each input data, we have 2 layers of CNN with one final fully connected layer. Convolution layer contains  8 and 16 feature maps for $X^k$ and $X^f$; and  4 and 8 feature maps for $X^{\overline{p}^{m}}$. We use dropout with rate 10 \% in each layer in order to prevent over-fitting.
Finally, we combine CNN output and perform two additional fully connected layers with size 200 and 1(one final output).
Presented algorithm is used to learn dependence between multi-input data and upscaled nonlinear transmissibilities.

For error calculation on the train and test dataset, we use  mean square errors, relative mean absolute and relative root mean square errors
\[
MSE = \sum_i |Y_i - \tilde{Y}_i|^2,
\quad
RMSE = \sqrt{ \frac{\sum_i |Y_i - \tilde{Y}_i|^2 }{\sum_i |Y_i|^2 } },
\quad
MAE = \frac{\sum_i |Y_i - \tilde{Y}_i|  }{\sum_i |Y_i|},
\]
where $Y_i$  and $\tilde{Y}_i$ denotes reference and predicted values for sample $X_i$

\begin{figure}[h!]
\centering
\includegraphics[width=0.7 \textwidth]{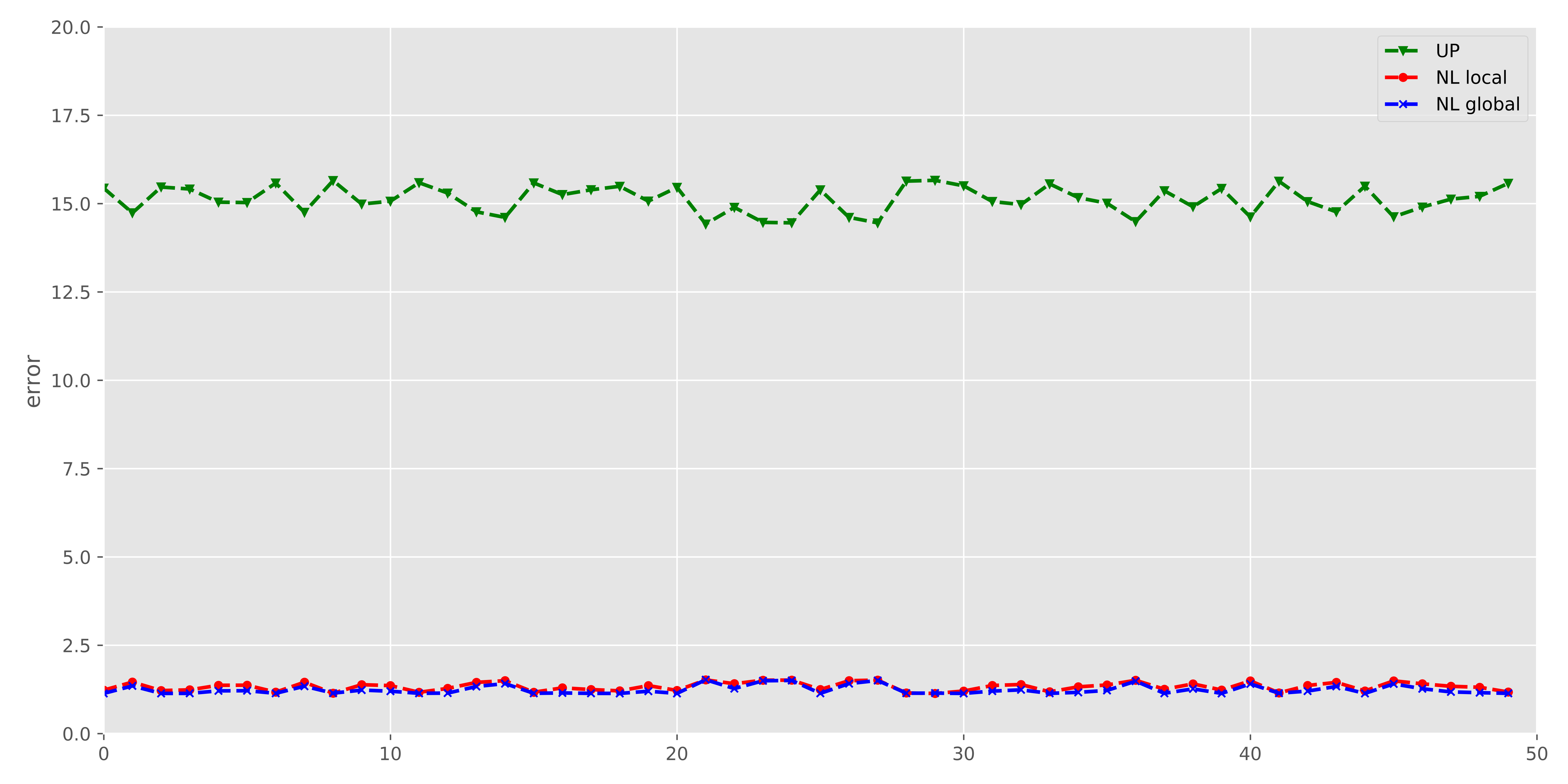}
\caption{Upscaling error for coarse grid parameters predicted using machine learning algorithm for \textit{Test 1}.
Green color: $e(\overline{u}^{UP})$. Red and blue colors:  $e(\overline{u}^{NL})$ with local and global calculations   }
\label{fig:err-ml}
\end{figure}

Convergence of the loss functions  for three neural networks for \textit{Test 1} are presented in Figure \ref{fig:ml1}, where we plot the MSE loss function vs epoch number for train and validation sets.
In Figure \ref{fig:ml2}, we present  a parity plots comparing reference values against predicted using trained neural networks for train and test datasets (green and blue colors).
Learning performance for neural networks are presented in Tables \ref{fig:err-ml} for global and local datasets. We observe good convergence of the relative errors  for train and test sets with $1-2 \%$ of RMSE.

Next, we consider errors between solution of the coarse grid problem with reference  and predicted upscaled transmissibilities.
In Figure \ref{fig:err-ml}, we present results for 50 test problems with random value of the source term. We show a relative $L_2$ errors for pressure head on the coarse mesh with classic upscaling algorithm and using new nonlocal nonlinear transmissibilities. We observe small errors ($1-2 \%$) for predicted nonlocal nonlinear transmissibilities compared with classical upscaling technique, where we have $\approx 15 \%$ of relative $L_2$ errors for pressure head.
Furthermore, we see that local calculation of the dataset provide similar results as a globally caclulated data.

\begin{figure}[h!]
\centering
\includegraphics[width=0.99\linewidth]{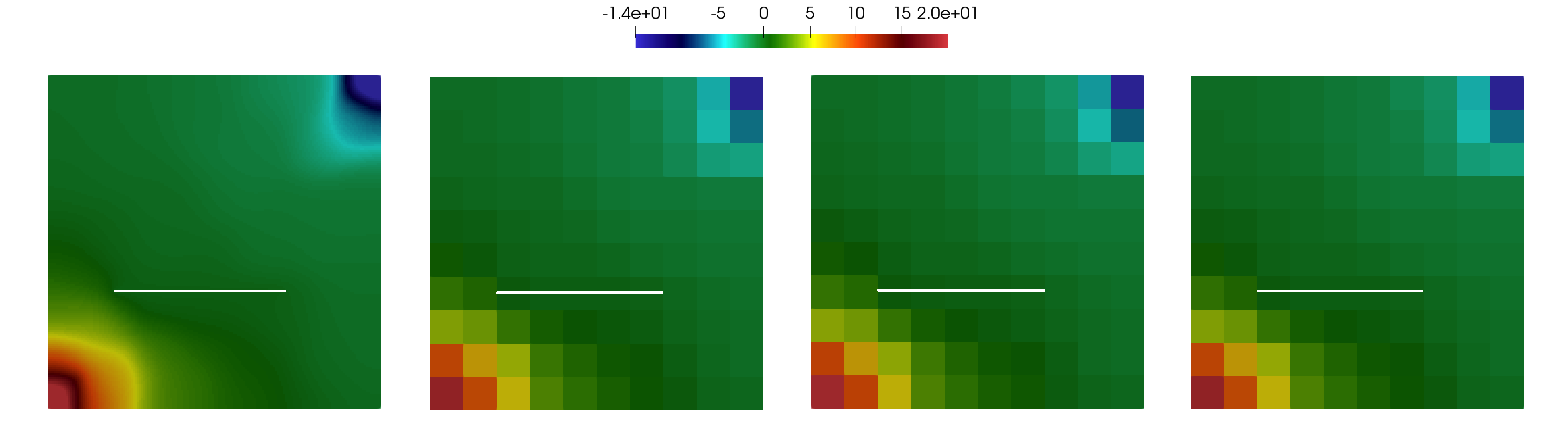}
\caption{Reference fine grid solution ($u^{fine}$), mean value on coarse grid of the fine grid solution ($\overline{u}^{fine}$), coarse grid solution using upscaling method ($\overline{u}^{UP}$)  and coarse grid solution using nonlinear nonlocal machine learning method ($\overline{u}^{NL}$). Nonlinear flow problem (\textit{Test 1}). Pressure on final time $t_m$, $m = 20$ }
\label{fig:uu}
\end{figure}

In Figure \ref{fig:uu}, we depict solution of the problem on the fine grid, coarse grid upscaled solution using classic approach from Section 3 and for new method presented in Section 4 ($u^{fine}$, $\overline{u}^{fine}$, $\overline{u}^{UP}$ and $\overline{u}^{NL}$).
For $\overline{u}^{UP}$, we apply presented upscaling method \ref{uns-Tup}, \ref{uns-Tup-fm} and \ref{uns-up1}.
We have $e(\overline{u}^{UP}) = 14.772 \%$  and $e(\overline{u}^{NL}) = 1.463 \%$  at final time.  For the nonlinear nonlocal transmissibilities, we set $\varepsilon = 0.5 \cdot 10^{-1}$ for $NN_1$ and  $NN_2$, $\varepsilon = 10^{-20}$ for $NN_3$. Note that, we didn't construct
$NN_4$ of data ($T^{ff, NL}$) because for our test problem we observe almost constant pressure on the fracture and set $T^{ff, NL} = T^{ff, UP}$ on the coarse grid.

We perform training of the neural networks on the GPU, where we train three neural networks: $NN_1$, $NN_2$ and $NN_3$.
Online stage (neural network training) time is 25 minutes for $NN_1$, 28 minutes for $NN_2$ and 6 minutes for $NN_3$ on GPU (GeForce GTX 1060). Note that, the training time depends on size of the dataset and GPU model. Here we didn't consider time of the dataset construction which depends on type of calculations (global or local) and  number of solution snapshots, that we used for training. Number of snapshots ($N_r$)  is also effects to the algorithm errors because we should have sufficient number of snapshots to capture all  variations of the input data to know how it effects to the output.

Time of the online stage contains 6.6 seconds of loading three neural networks and 13.0 seconds for calculations on the $10 \times 10$ coarse grid with prediction of the nonlinear nonlocal transmissibilities.
Fine grid calculations time is 454 seconds for 20 time steps on $640 \times 640$ fine grid.
We have approximately 35 time faster calculations for a new method with small error of the coarse grid solution.


\subsection{Nonlinear flow and transport problem }

We consider solution of the two-phase flow problem in fractured and heterogeneous porous media. For nonlinear coefficient, we set
$\lambda^w(s)  = s^2$ and $\lambda^n(s)  = (1-s)^2$.
In Figure \ref{fig:kx} (third column),  we show the heterogeneous porous matrix permeability $k^m(x)$ and fracture position. We set $\phi^{\alpha} = 1$ ($\alpha=m,f$), $k^f = 10^3$ and $T_{max} = 25 \cdot 10^{-3}$ with 250 time steps. Coarse grid is $10 \times 10$ and fine grid is $160 \times 160$ for domain $\Omega$.

\begin{figure}[h!]
\centering
\includegraphics[width=0.48\linewidth]{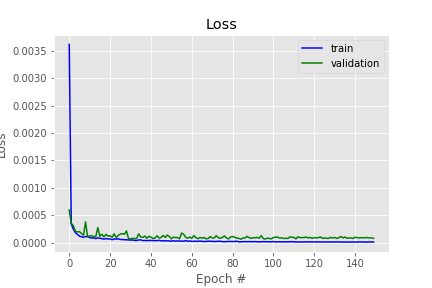}
\includegraphics[width=0.48\linewidth]{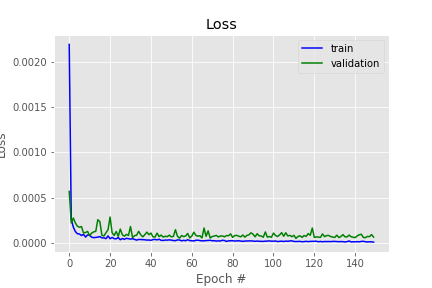}\\
\includegraphics[width=0.48\linewidth]{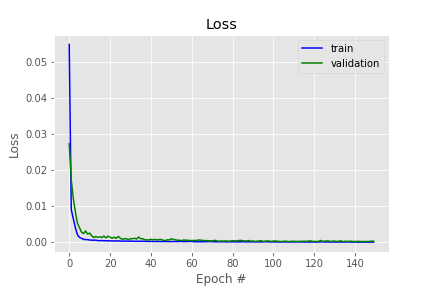}
\includegraphics[width=0.48\linewidth]{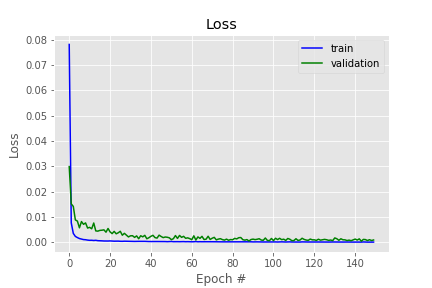}
\caption{Learning process  for \textit{Test 2} (nonlinear flow and transport). Loss function vs epoch.
First row: $NN_1$ and  $NN_2$ for vertical  and horizontal ${T}^{mm, NL}$.
Second row: $NN_3$  and $NN_4$ for  ${T}^{mf, NL}$  and ${T}^{ff, NL}$}
\label{fig:ml-t2}
\end{figure}

\begin{table}[h!]
\begin{center}
\begin{tabular}{ | c | c c c | }
\hline
& MSE & RMSE  (\%)  & MAE (\%) \\
\hline
$NN_1$ 	& 0.017 & 1.316 & 0.959 \\
$NN_2$ 	& 0.043 & 2.092 & 1.507 \\
$NN_3$ 	& 0.014 & 1.218 & 0.778\\
$NN_4$ 	& 0.052 & 2.301 & 1.328 \\
\hline
\end{tabular}
\end{center}
\caption{Learning performance of machine learning algorithm  for \textit{Test 2} (nonlinear flow and transport). Errors for train and test sets }
\label{tab:ml-t2}
\end{table}

For the training of the neural networks, we use a global dataset, where we extract local information from the fine grid calculations on the global domain $\Omega$. For generation of the train datasets, we use a three random shapshots ($N_r = 3$) with $T_{max} = 40 \cdot 10^{-3}$ and 400 time steps.
We train four neural networks for each type of transmissibility: $NN_1$ for horizontal coarse edges for matrix-matrix flow, $NN_2$ for vertical coarse edges s for matrix-matrix flow, $NN_3$ for matrix - fracture flow and $NN_4$ for fracture - fracture flow.
The train dataset for first and second neural networks contains $N = 108 000$; $N = 6 000$ for $NN_3$ and $N = 4 800$ for $NN_4$, where dataset is randomly divided into training and validation sets with $80:20$ ratio.
Each sample $X_l$ contains information about heterogeneous permeability and fracture position up to fine grid resolution in local domain, mean value of the solution in oversampled local domain (coarse grid)
\[
X_l = (X_l^k, X_l^f, X_{l+}^{\overline{p}^{\alpha}}, X_{l+}^{\overline{s}^{\alpha}}, X_{l+}^{\overline{p}^{\beta}}, X_{l+}^{\overline{s}^{\beta}})
\]
and output
\[
Y_l = (T_l^{\alpha \beta, NL}, T_l^{w, \alpha \beta, NL}), \quad \alpha,\beta = m,f.
\]

For calculations, we use 150 epochs with a batch size $N_b = 90$ and perform calculations on GPU. Architecture of the neural networks are similar to the previous test problem but as output for this case, we obtain two values, $T$.
Learning performance for neural networks are presented in Tables \ref{fig:ml-t2} and \ref{tab:ml-t2} for train datasets. We observe a good convergence with small error for each neural network.

\begin{figure}[h!]
\centering
\includegraphics[width=0.48 \textwidth]{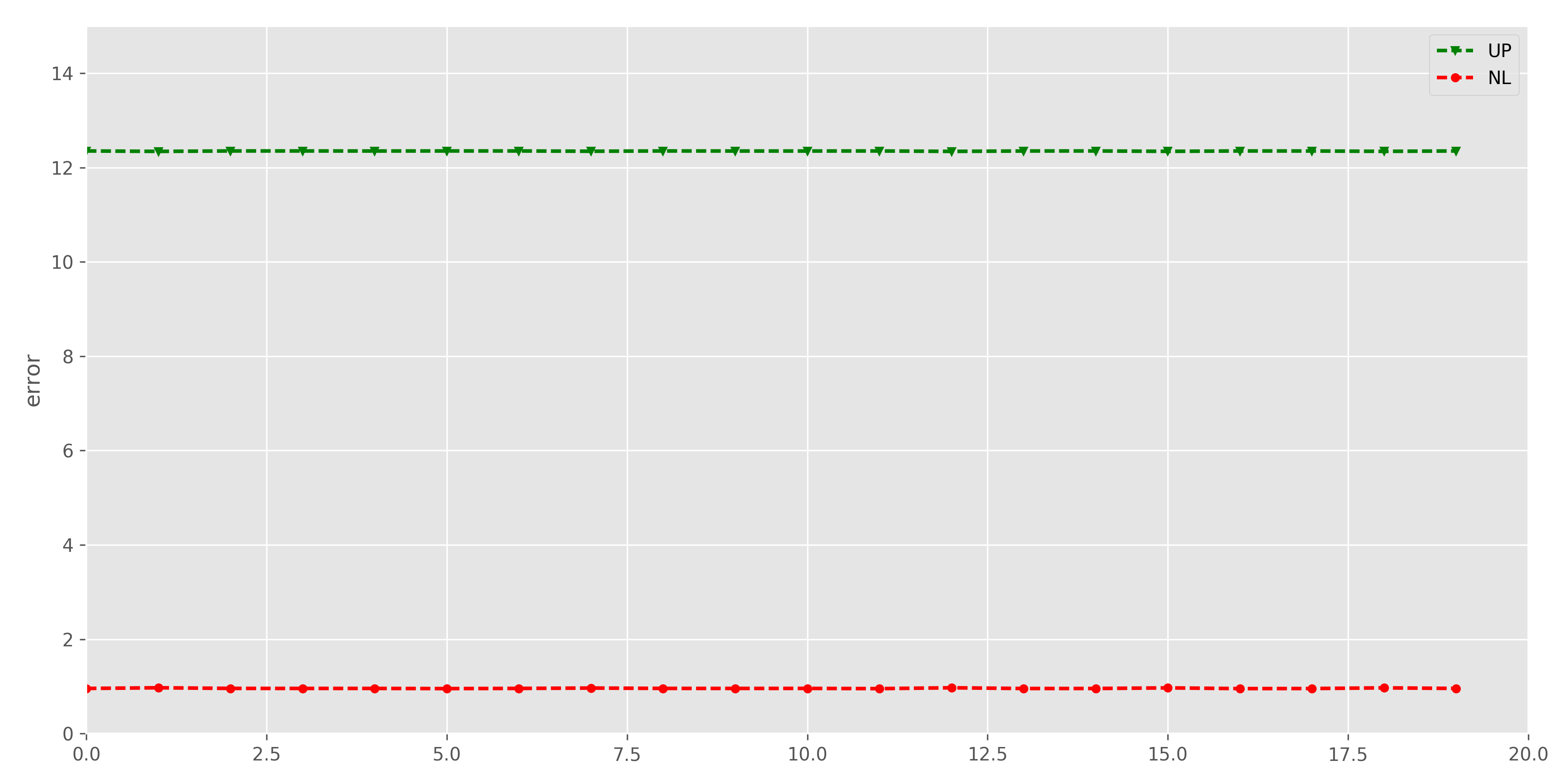}
\includegraphics[width=0.48 \textwidth]{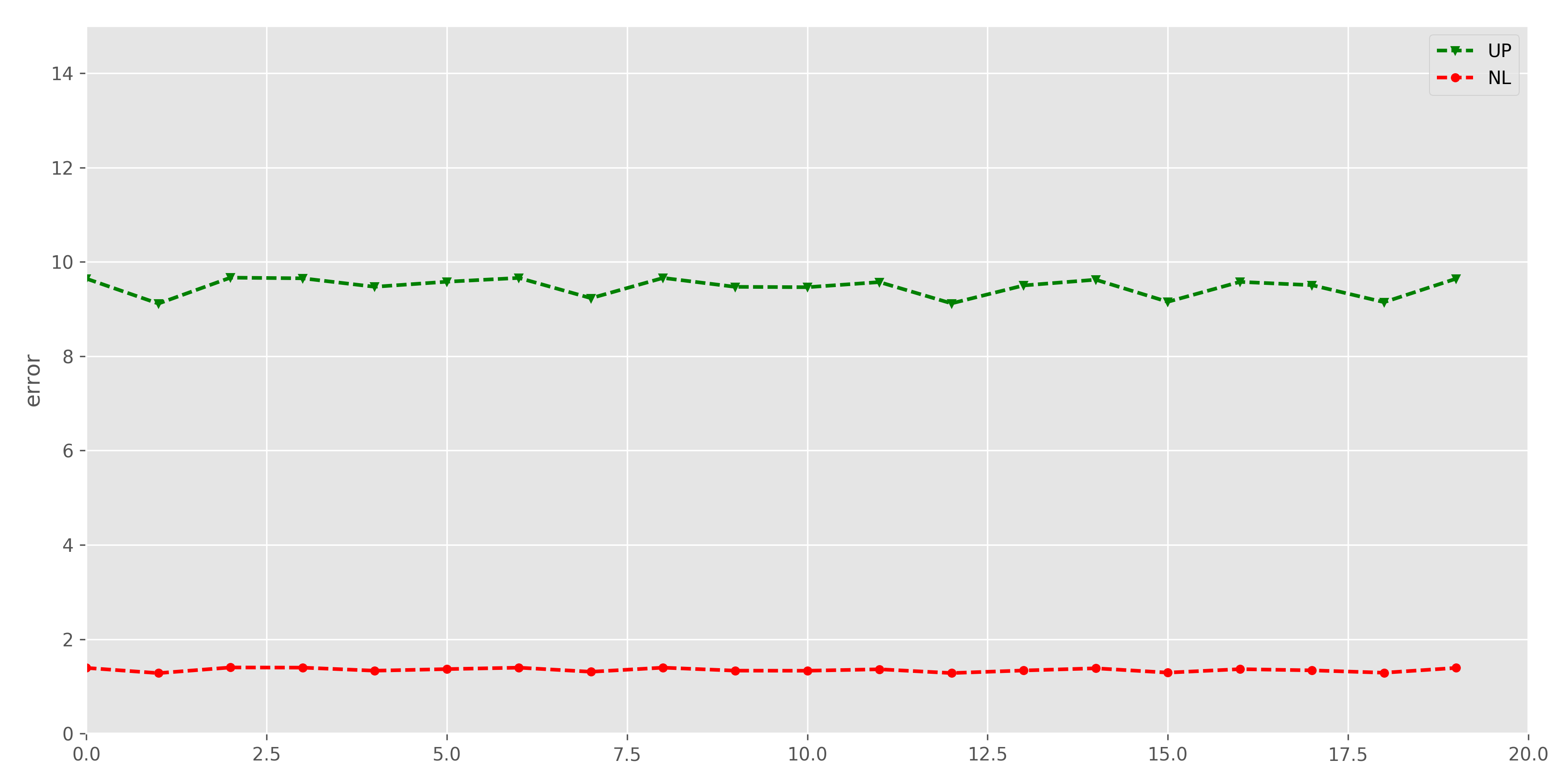}\\
\includegraphics[width=0.48 \textwidth]{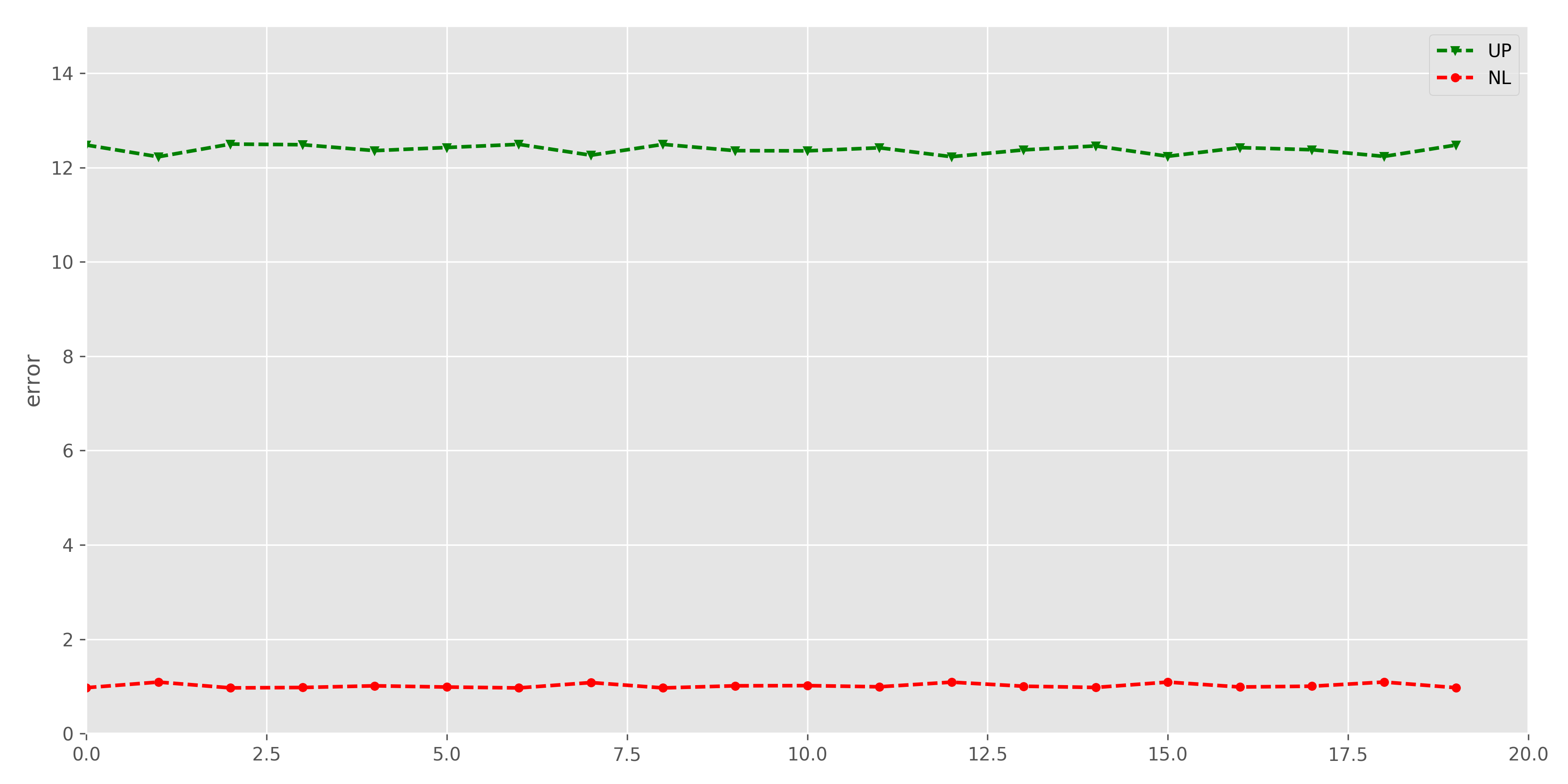}
\includegraphics[width=0.48 \textwidth]{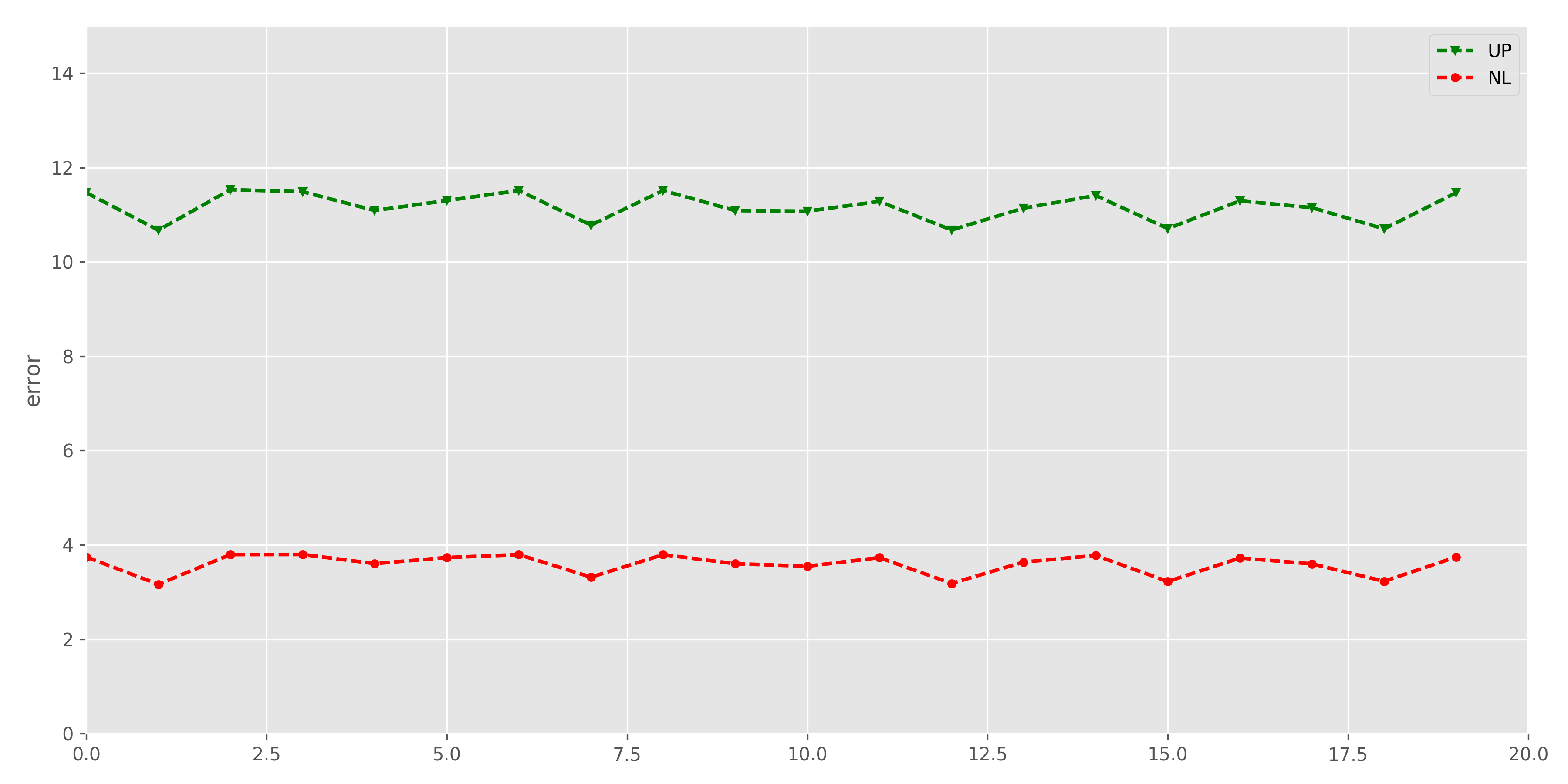}\\\includegraphics[width=0.48 \textwidth]{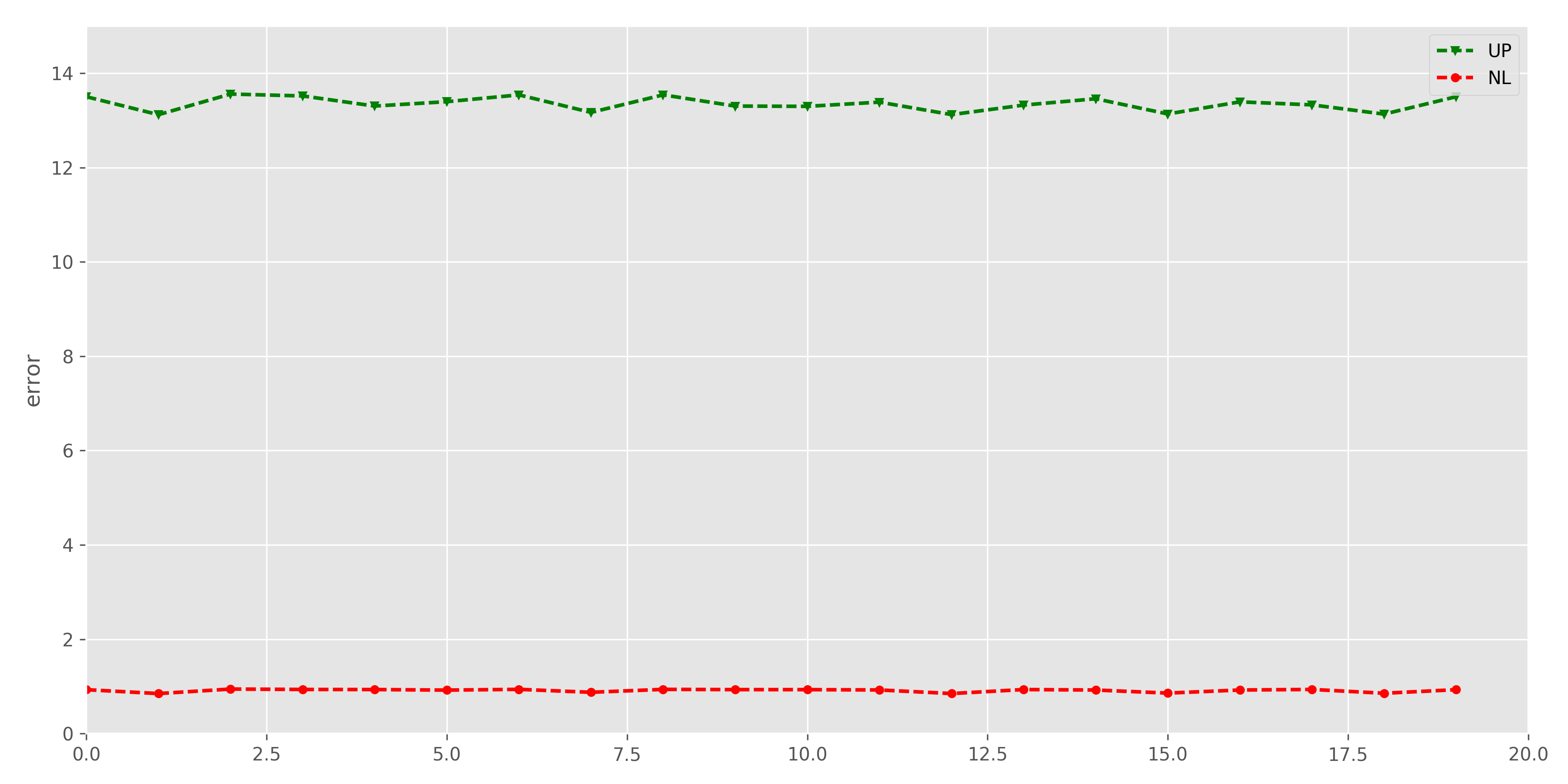}
\includegraphics[width=0.48 \textwidth]{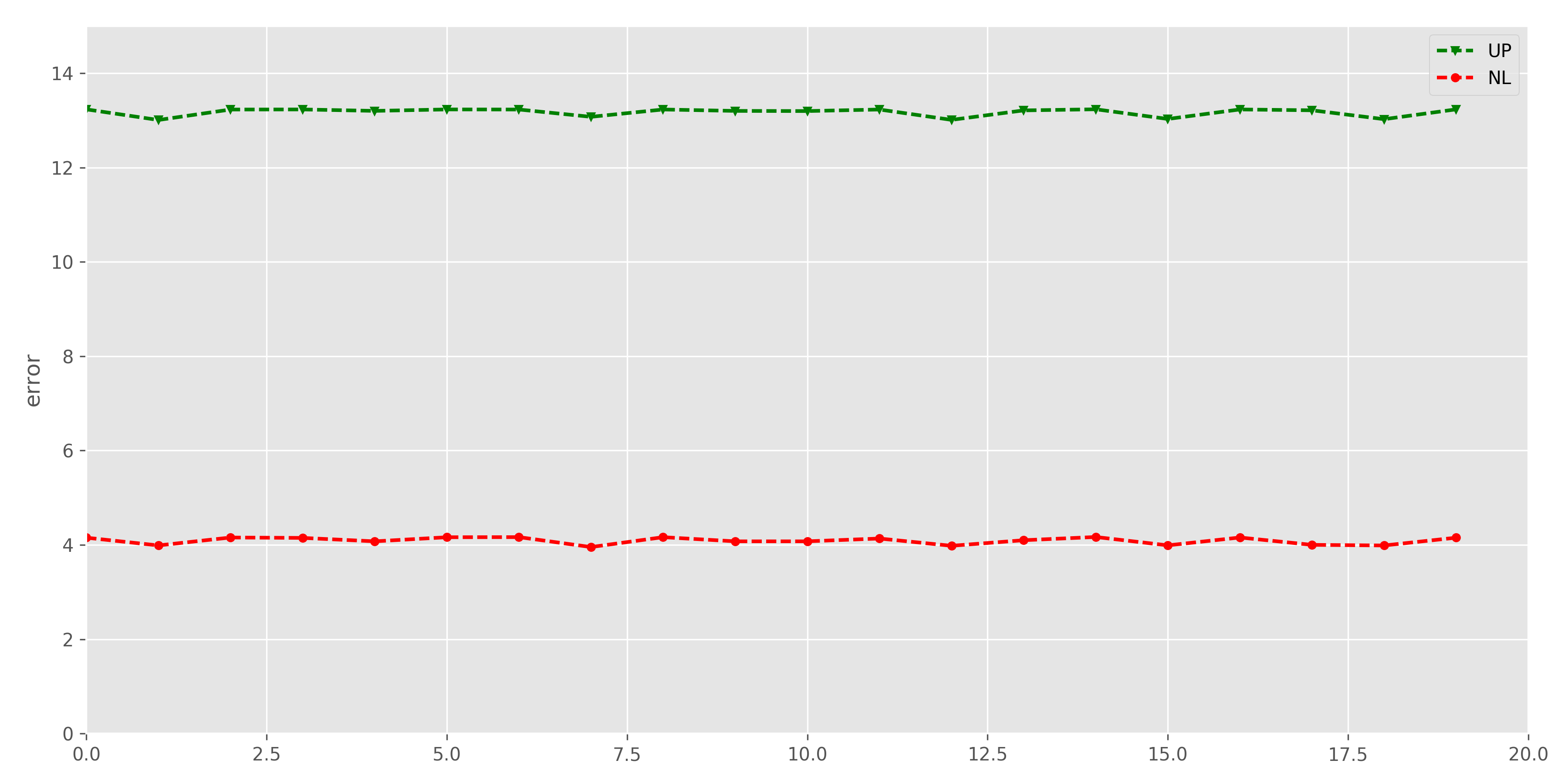}
\caption{Upscaling error for coarse grid parameters predicted using machine learning algorithm  for \textit{Test 2} (nonlinear flow and transport).
Green color: $e(\overline{u}^{UP})$. Red color:  $e(\overline{u}^{NL})$ }
\label{fig:err-ml-t2}
\end{figure}

In Figure \ref{fig:err-ml-t2}, we present results for 20 test problems with random value of the source terms. We show a relative mean square error in percentages for pressure and for saturation on the coarse mesh with classic upscaling algorithm and using new nonlocal nonlinear transmissibilities.

\begin{figure}[h!]
\centering
\includegraphics[width=0.99\linewidth]{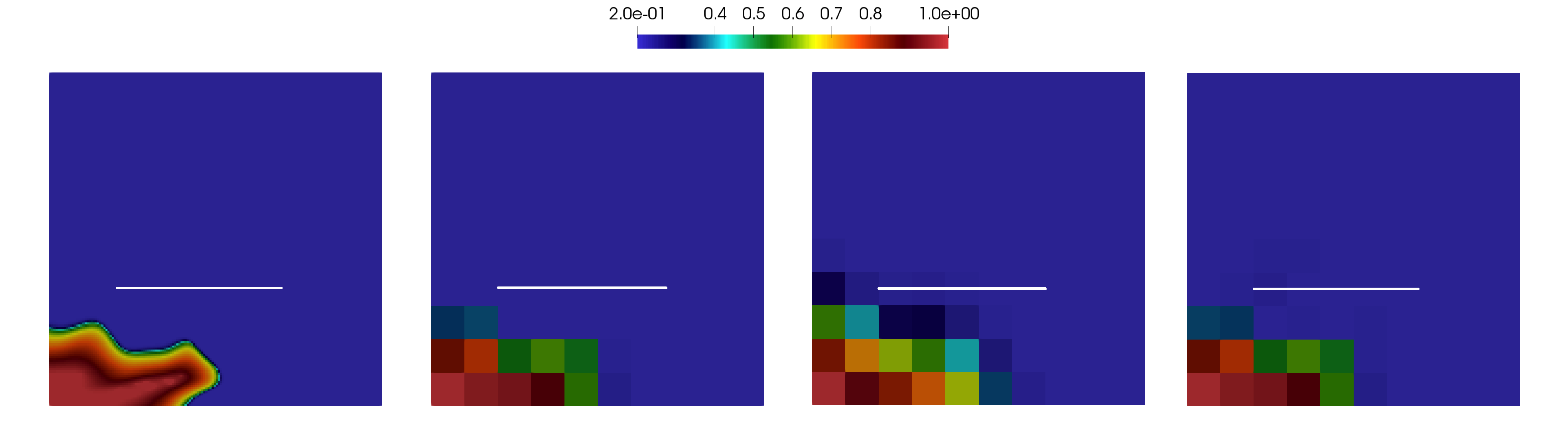}\\
\includegraphics[width=0.99\linewidth]{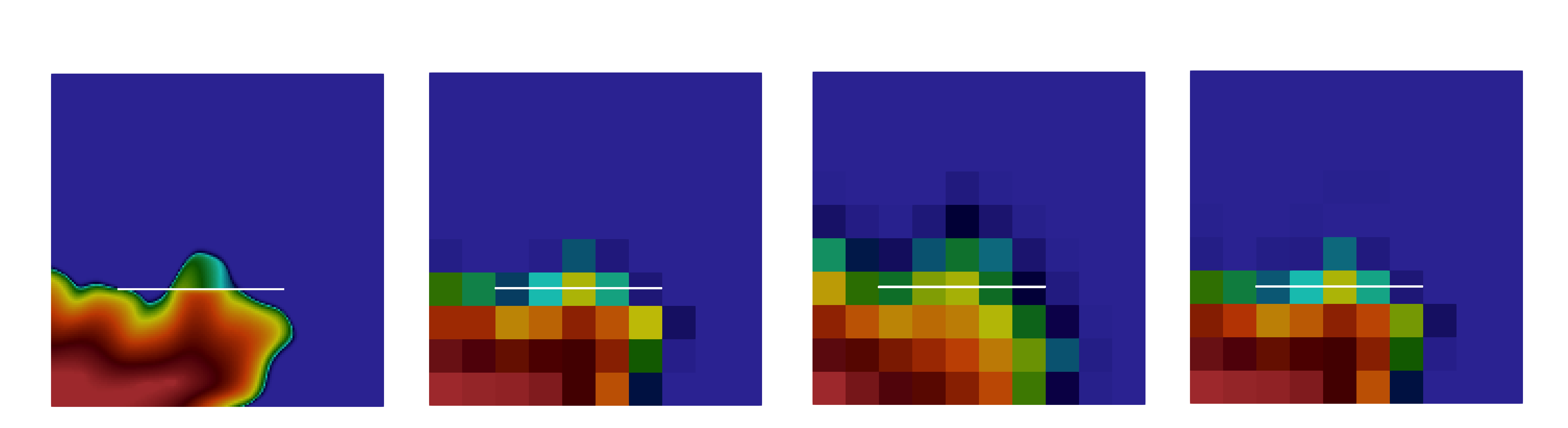}\\
\includegraphics[width=0.99\linewidth]{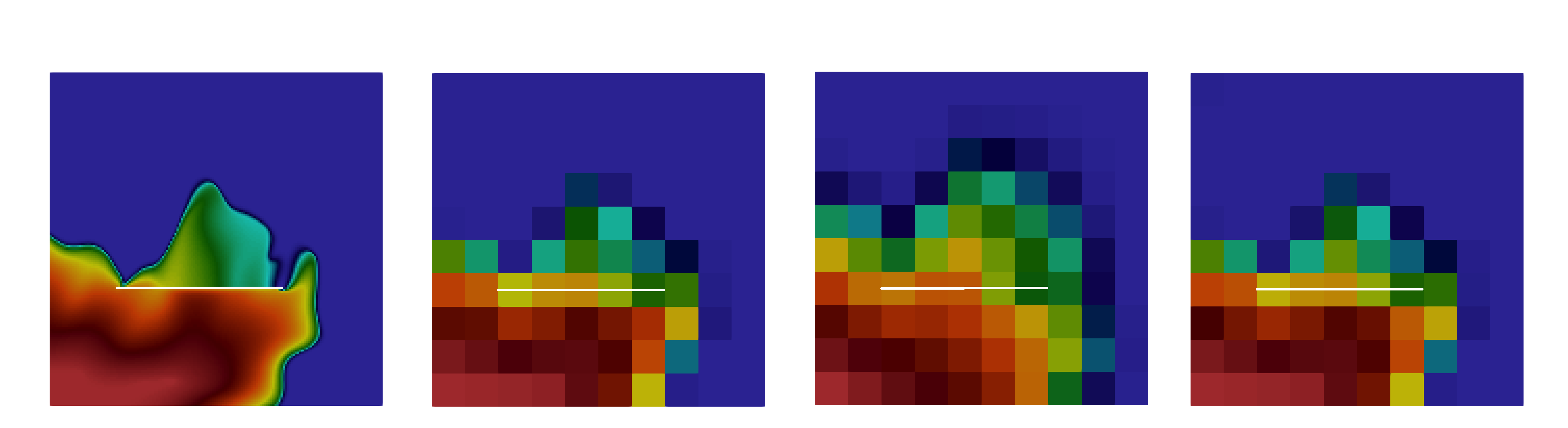}\\
\includegraphics[width=0.99\linewidth]{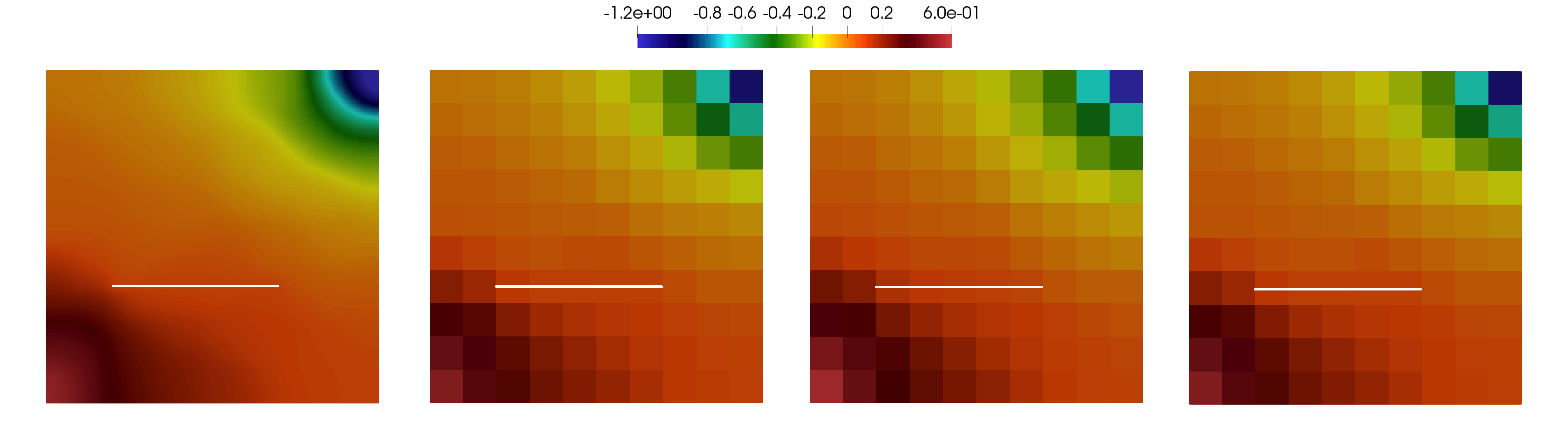}
\caption{Reference fine grid solution ($s^{fine}$, $p^{fine}$), mean value on coarse grid of the fine grid solution ($\overline{s}^{fine}$, $\overline{p}^{fine}$), coarse grid solution using upscaling method ($\overline{s}^{UP}$, $\overline{p}^{UP}$)  and coarse grid solution using nonlinear nonlocal machine learning method ($\overline{s}^{NL}$, $\overline{p}^{NL}$).
Nonlinear flow and transport problem (\textit{Test 2}).
First row: saturation for time $t_m$, $m = 50$.
Second row: saturation for time $t_m$, $m = 150$.
Third row: saturation for time $t_m$, $m = 250$.
Fourth row: pressure for time $t_m$, $m = 250$ }
\label{fig:uu-t2}
\end{figure}

In Figure \ref{fig:uu-t2}, we depict solution of the problem using different methods.
On the first conlumn, we depict a reference fine grid solution ($s^{fine}$, $p^{fine}$), mean value on coarse grid of the fine grid solution ($\overline{s}^{fine}$, $\overline{p}^{fine}$) on the second column, coarse grid solution using upscaling method ($\overline{s}^{UP}$, $\overline{p}^{UP}$) on the third column and coarse grid solution using nonlinear nonlocal machine learning method ($\overline{s}^{NL}$, $\overline{p}^{NL}$) on the fourth column.
On the first, second and third rows, we show a saturation for time $t_m$, $m = 50, 150, 250$ and on fourth row, we have pressure for time $t_m$, $m = 250$.
For solution on the coarse grid ($\overline{p}^{UP}$ and $\overline{s}^{UP}$), we applied classic upscaling method (see Section 3).
Fine grid (reference) solution is performed using finite volume approximation with embedded discrete fracture model, where for error calculations we used a mean values of the reference solution on the coarse grid, $\overline{p}^{fine}$ and $\overline{s}^{fine}$.
On the last column of the Figure \ref{fig:uu-t2}, we depict a coarse grid solution using nonlinear nonlocal transmissibilities that calculate based on the machine learning approach.
For machine learning approach, we have
$e(\overline{p}^{NL}) = 0.920 \%$, $ e(\overline{s}^{NL}) = 3.957 \%$, and for upscaling
$e(\overline{p}^{UP}) = 13.519 \%$, $ e(\overline{s}^{UP}) = 13.227 \%$ at final time $t_m$, $m = 250$.
For the nonlinear nonlocal transmissibilities, we set $\varepsilon_S = 10^{-2}$ for transport and  $\varepsilon = 0.5 \cdot 10^{-2}$ for $NN_1$, $\varepsilon = 10^{-4}$ for $NN_2$, $\varepsilon = 10^{-3}$ for $NN_3$ and $\varepsilon = 10^{-20}$ for $NN_4$ for flow.

We perform training of the neural networks on the GPU, where we train four neural networks: $NN_1$, $NN_2$, $NN_3$ and $NN_4$. Online stage (neural network training) time is 80 minutes for $NN_1$, 59 minutes for $NN_2$, 2 minutes for $NN_3$ and 4 minutes for $NN_4$ on GPU (GeForce GTX 1060). Note that, the training time depends on size of the dataset and GPU model.
Time of the online stage contains 16.7 seconds of loading four neural networks and 46.9 seconds for calculations on the $10 \times 10$ coarse grid with prediction of the nonlinear nonlocal transmissibilities.
Fine grid calculations time is 812 seconds for 250 time steps on $160 \times 160$ fine grid for transport and flow model.
We observe a good results with fast calculations using a machine learning algorithm for presented method.



\section{Conclusion}

In this work, we consider two nonlinear problems in heterogeneous and fractured porous media. Mathematical models are formulated as a general multicontinuum models, where fine grid approximations are constructed using finite volume method.
For the accurate solution of the nonlinear problems on the coarse grid, a novel machine learning algorithm combined with nonlinear nonlocal multicontinua approach for calculating nonlocal nonlinear transmissibilities is presented and investigated. We presented the construction of the dataset for training deep neural networks. The construction of the neural network is based on the multi-input convolutional neural networks, where GPU is used for performing a fast learning process.
To illustrate the applicability of the presented method, we presented numerical results for two test problems. Numerical results showed that presented algorithm provides fast and accurate calculations of the nonlocal nonlinear transmissibilities.

\section{Acknowledgements}
MV's  work is supported by the mega-grant of the Russian Federation Government N14.Y26.31.0013 and RSF N17-71-20055.
EC's work is partially supported by the Hong Kong RGC General Research Fund (Project numbers 14304217 and 14302018) and the CUHK Faculty of Science Direct Grant 2018-19.

\bibliographystyle{plain}
\bibliography{lit}

\end{document}